\documentstyle[amscd]{amsart}

\newtheorem{theorem}{Theorem}[section]
\newtheorem{lemma}[theorem]{Lemma}
\newtheorem{proposition}[theorem]{Proposition}
\newtheorem{corollary}[theorem]{Corollary}

\theoremstyle{definition}

\theoremstyle{remark}
\newtheorem{ack}{Acknowledgment}  

\newcommand{\Ext}{\operatorname{Ext}}
\newcommand{\Hom}{\operatorname{Hom}}
\newcommand{\ob}{\operatorname{ob}}

\newcommand{\OMod}{\operatorname{Mod_{\cal O}}}
\newcommand{\grsmod}{\operatorname{Mod}_{\mit{S},\operatorname{gr}}}
\newcommand{\grssmod}
	{\operatorname{Mod}_{\mit{S}\otimes\mit{S},\operatorname{gr}}}
\newcommand{\coker}{\operatorname{coker}}
\newcommand{\im}{\operatorname{im}}
\newcommand{\Id}{\operatorname{Id}}

\newcommand{\tr}{\operatorname{tr}}
\newcommand{\Spec}{\operatorname{Spec}}

\newcommand{\length}{\operatorname{length}}
\newcommand{\sat}{\operatorname{sat}}
\newcommand{\End}{\operatorname{End}}

\newcommand{\Pnt}{\Bbb{P}^{n+3}}
\newcommand{\PN}{\Bbb{P}^N}

\newsymbol\twoheadrightarrow 1310
\newsymbol\boxtimes 1202

\begin{document}

\title{Pfaffian Subschemes}
\author{Charles H. Walter}
\address{URA 168\\Math\'ematiques\\Universit\'e de
	Nice--Sophia-Antipolis\\Parc Valrose\\F-06108 NICE Cedex
02\\France}
\email{walter@@math.unice.fr}
\thanks{The author was supported in part by NSA research grant
	MDA904-92-H-3009.}

\maketitle

\begin{abstract}
A subscheme $X\subset \Pnt$ of codimension $3$ is {\em Pfaffian} if it
is the degeneracy locus of a skew-symmetric map $f:\cal{E}\spcheck(-t)
@>>> \cal{E}$ with $\cal{E}$ a locally free sheaf of odd rank on
$\Pnt$.  It is shown that a codimension $3$ subscheme $X\subset\Pnt$
is Pfaffian if and only if it is locally Gorenstein, subcanonical
(i.e.\ $\omega_X\cong\cal O_X(l)$ for some integer $l$), and the
following parity condition holds: if $n\equiv 0 \pmod{4}$ and $l$ is
even, then $\chi (\cal O_X (l/2))$ is also even.

The paper includes a modern version of the Horrocks correspondence,
stated in the language of derived categories.  A local analogue of the
main theorem is also proved.

\end{abstract}

One method of constructing a codimension $3$ subscheme $X \subset
\Pnt$ is to consider a {\em Pfaffian subscheme}. (This is discussed in
a recent paper of Okonek \cite{O}.)  That is, one considers the
degeneracy locus of a skew-symmetric map $f:\cal{E}\spcheck(-t) @>>>
\cal{E}$ such that $\cal{E}$ is a locally free sheaf of odd rank
$2p+1$ on $\Pnt$, and $f$ is generically of rank $2p$ and degenerates
to rank $2p-2$ in the expected codimension $3$.  The scheme $X$ then
has a locally free resolution of the form
\begin{equation}
\label{resol}
0 @>>> \cal{O}_{\Pnt}(-t-2s) @>h>> \cal{E}\spcheck(-t-s)
@>f>> \cal{E}(-s) @>g>> \cal{O}_{\Pnt} @>>> \cal O_X
\end{equation}
where $s=c_1(\cal{E})+pt$, and where $g$ and $h=g\spcheck(-t-2s)$ are
given locally by the Pfaffians of order $2p$ of $f$.  This resolution
is just a patching together of the local version studied in \cite{BE}.
The self-duality of the resolution (\ref{resol}) implies that $X$ is
locally Gorenstein with canonical sheaf $\omega_X\cong\cal O_X
(t+2s-n-4)$.  Thus Pfaffian subschemes are always locally Gorenstein
of codimension $3$ in $\Pnt$ and are subcanonical, i.e.\ they satisfy
$\omega_X\cong\cal O_X (l)$ for some integer $l$.

It is now natural to consider the following question asked by Okonek
\cite{O}:  Are all locally Gorenstein subcanonical subschemes of
codimension $3$ in $\Pnt$ Pfaffian?  Arithmetically Gorenstein
subschemes of codimension $3$ certainly are Pfaffian because of a
structure theorem of Buchsbaum and Eisenbud (\cite{BE} Theorem 2.1).
We will show that in conjunction with a certain number of other ideas,
their method can be adapted to yield the following result:

\begin{theorem}
\label{main}
Let $k$ be a field not of characteristic $2$.  Suppose $X\subset
\Pnt_k$ is a locally Gorenstein subscheme of equidimension $n>0$ such
that $\omega_X\cong \cal O_X (l)$ for some integer $l$.  Then $X$ is a
Pfaffian subscheme if and only if the following parity condition
holds\rom: if $n\equiv 0 \pmod{4}$ and $l$ is even, then $\chi (\cal
O_X (l/2))$ is also even.
\end{theorem}

As Okonek pointed out, the Barth-Lefschetz theorems imply that every
smooth subvariety of codimension $3$ in $\Bbb{P}^8$ and $\Bbb{P}^9$ is
subcanonical.  Hence we have a corollary:

\begin{corollary}
Every smooth subvariety of codimension $3$ in $\Bbb{P}^8$ and
$\Bbb{P}^9$ is a Pfaffian subscheme.
\end{corollary}

Of course the corollary holds for $\PN$ with $N\geq 10$ as well (with
certain conditions if $N\equiv 3 \pmod{4}$). But in this range all
smooth subvarieties of codimension $3$ in $\PN$ are supposed to be
complete intersections according to Hartshorne's Conjecture.

The parity condition in Theorem \ref{main} may be deduced as follows.
Suppose $X$ is Pfaffian and $n$ and $l$ are both even.  We may twist
the resolution (\ref{resol}) by $l/2$ to get
\[
0 @>>> \omega_{\Bbb{P}}(-l/2) @>>> \cal{E}\spcheck \otimes \omega
_{\Bbb{P}}(s-l/2) @>>> \cal{E}(l/2-s) @>>> \cal{O}_{\Bbb{P}}(l/2) @>>>
\cal O_X (l/2) @>>> 0.
\]
{}From this sequence and Serre duality on $\Pnt$, it now follows that:
\[
\chi (\cal O_X (l/2)) = 2\chi (\cal{O}_{\Bbb{P}}(l/2)) -2\chi
(\cal{E}(l/2-s)) \equiv 0\pmod 2.
\]

The main theorem has the following local analogue.  Call an unmixed
ideal $I$ of height $3$ in a regular local ring $(R,\frak m,k)$ {\em
Pfaffian} if $R/I$ has a resolution of the form
\[
0 @>>> R @>{g\spcheck}>> E\spcheck @>f>> E @>g>> R @>>> R/I
\]
with $E$ is a reflexive $R$-module of odd rank such that $E_{\frak p}$
is a free $R_{\frak p}$-module for all prime ideals $\frak p\neq \frak
m$, and $f$ skew-symmetric.  The canonical module of any unmixed ideal
$I$ of height $3$ is by definition $\omega_{R/I} = \Ext^3_R(R/I,R)$.
It is always saturated.  We prove the following:

\begin{theorem}
\label{RLR}
Let $(R,\frak m,k)$ be a regular local ring of dimension $n>4$ with
residue field not of characteristic $2$.  Let $I$ be an unmixed ideal
of $R$ of height $3$.  Then $I$ is Pfaffian if and only if the
following three conditions hold\rom:

\rom(a\rom) $(R/I)_{\frak p}$ is
Gorenstein for all prime ideals $\frak p\neq\frak m$,

\rom(b\rom)
$\omega_{R/I}\cong (R/I)^{\sat}$, and

\rom(c\rom) if $n\equiv 0\pmod
4$, then $H^{n/2}_{\frak m}(I)$ is of even length.
\end{theorem}

In both theorems we have stated the parity condition only for $n\equiv
0\pmod{4}$ rather than for all even $n$.  This is because the graded
commutativity built into cohomology rings causes the perfect pairing
of Serre duality
\begin{equation}
\label{pairing}
H^{n/2}(\cal O_X (l/2)) \times H^{n/2}(\cal O_X (l/2)) @>>> H^n(\cal
O_X (l))\cong k
\end{equation}
(or its analogue in local duality) to be $(-1)^{n/2}$-symmetric. Hence
if $n\equiv 2\pmod{4}$, then $H^{n/2}(\cal O_X (l/2))$ admits a
non-degenerate skew-symmetric bilinear form and so is of even
dimension.  Thus subcanonical varieties with $n\equiv 2 \pmod{4}$ and
$l$ even automatically have
\[
\chi (\cal O_X (l/2)) \equiv h^{n/2}(\cal O_X (l/2))\equiv 0 \pmod{2}.
\]

In characteristic $2$ the perfect pairing (\ref{pairing}) and its
local analogue apparently need not be alternating even if $n\equiv
2\pmod 4$.  Thus one cannot expect Theorem \ref{main} or
\ref{RLR} to be valid in characteristic $2$ unless the phrase ``if
$n\equiv 0\pmod 4$'' is replaced by the phrase ``if $n$ is even.''
However, we will show that with this modification, both theorems are
valid in characteristic $2$.

\subsection{Outline of the Paper}

In the first section we review the proof of the local version of
Theorem \ref{main} given by Buchsbaum and Eisenbud (\cite{BE} Theorem
2.1).  We show that their proof will work for us if we can replace
their minimal projective resolution by a locally free resolution of
$\cal O_X$ which satisfies two properties (Proposition
\ref{conditions}).  The rest of the paper is devoted to finding a
locally free resolution of $\cal O_X$ which satisfies these
properties.

Our main tool for constructing this locally free resolution is the
Horrocks correspondence of \cite{H}.  In the second section of the
paper, we give a modern description of this correspondence using
derived categories.  This point of view is not identical to Horrocks',
so we have felt it prudent to include a full proof of Horrocks'
principal result (Theorem \ref{Horrocks}) from this point of view.
However, the derived categories viewpoint is useful because it permits
us to further develop Horrocks' ideas so as to obtain a method for
transfering a portion of the cohomology of the coherent sheaf $\cal
O_X$ to a locally free sheaf in a controlled way (Proposition
\ref{functorial}).  This is critical for our construction.

In the third section we apply the Horrocks correspondence to construct
a particular locally free resolution of the form (\ref{resol}).  The
basic idea is to cut in half the cohomology of the subscheme $X$ by
using truncations of $\bold R\Gamma_*(\cal I_X)$.  Our results on the
Horrocks correspondence then permit us to find a vector bundle $\cal
F_1$ whose intermediate cohomology is one of the halves of the
cohomology of $\cal O_X$.  Moreover, there is a natural morphism from
this $\cal F_1$ to $\cal I_X$.  This more or less gives the right half
of the resolution, and the left half comes from the conventional
methods of the Serre correspondence.  We then show that if the
cohomology of $\cal O_X$ was cut in half properly (viz.\ if the
subcomplex carries an ``isotropic'' half of the cohomology), then the
resolution is self-dual in a very strong way: i.e.\ any chain map from
the resolution to its dual which extends the identity on $\cal{O}_X$
is necessarily an isomorphism of complexes.  This is one of the
properties required of the locally free resolution in order to make
the Buchsbaum-Eisenbud proof work.

In the fourth section we show that our locally free resolution of
$\cal O_X$ can be endowed with a commutative differential graded
algebra structure.  This is a matter of calculating the obstruction to
the lifting of a certain map.  This is the second property required of
the locally free resolution in order for the Buchsbaum-Eisenbud proof
to work.  This will complete the proof of Theorem \ref{main}.

In the fifth section we consider Theorem \ref{main} in characteristic
$2$.  Essentially, certain lemmas in the fourth section fail in
characteristic $2$ and must be replaced by analogues which are
slightly different.

In the sixth section we consider the results for regular local rings.
Theorem \ref{main} concerning projective spaces has an obvious
analogue (Theorem \ref{punc:spec}) for the punctured spectrum of a
regular local ring.  We show that this analogue is equivalent to
Theorem \ref{RLR}.

\begin{ack}
The author would like to thank R.~M.~Mir\`o-Roig who brought the
problem to his attention and with whom he had several discussions
concerning it.  The paper was written in the context of the Space
Curves group of Europroj.
\end{ack}

\section{The Buchsbaum-Eisenbud Proof}

In this section we review Buchsbaum and Eisenbud's proof of the local
version of Theorem \ref{main}. In particular, we describe the two
conditions that a locally free resolution of $\cal O_X$ must satisfy
in order for their proof to show that a subcanonical subscheme
$X\subset \Pnt$ is Pfaffian (Proposition \ref{conditions}).

\begin{theorem}[\cite{BE} Theorem 2.1]
\label{local}
Let $R$ be a regular local ring and $I$ an ideal of $R$ of height $3$
such that $R/I$ is a Gorenstein ring.  Then $I$ has a minimal
projective resolution of the form
\[
0 @>>> R @>{g\spcheck}>> F\spcheck @>f>> F @>g>> R @>>> R/I
\]
such that $F$ of odd rank $2p+1$, the map $f$ is skew-symmetric, and
$g$ is composed of the Pfaffians of order $2p$ of $f$.
\end{theorem}

\begin{pf*}{Sketch of Buchsbaum and Eisenbud's proof of Theorem
\ref{local}}
One considers a minimal projective resolution of $R/I$.  Since $R/I$
is Gorenstein, it is of the form
\[
\bold{P}^*: \qquad\quad 0 @>>> R @>{d_3}>> F_2 @>{d_2}>> F_1 @>{d_1}>>
R
\]
We now seek to find a way of identifying $F_2\cong F_1\spcheck$ so
that $d_2$ becomes skew-symmetric.

The first step is to endow $\bold{P}^*$ with the structure of a
commutative associative differential graded algebra (\cite{BE}
pp.~451--453).  To define the multiplication, they define
$S_2(\bold{P}^*) = (\bold{P}^* \otimes \bold{P}^*)/M^*$ where $M^*$ is
the graded submodule of $\bold{P}^*\otimes\bold{P}^*$ generated by
\[
\{a\otimes b-(-1)^{(\deg a)(\deg b)}b\otimes a \mid a,b \text{
homogeneous elements of }\bold{P}^*\}.
\]
Using universal properties of projective modules, they then construct
a map of complexes $\Phi : S_2(\bold{P}^*) @>>> \bold{P}^*$ which
extends the multiplication $R/I \otimes R/I$ and which is the identity
on the subcomplex $R\otimes\bold{P}^* \subset S_2(\bold{P}^*)$.  This
makes $\bold{P}^*$ into a commutative differential graded algebra.
The associativity of this algebra follows from the fact that it of
length $3$, i.e. $P_n = 0$ for $n\geq 4$.

The next step (p.~455) is to note that the multiplication $F_i \otimes
F_{3-i} @>>> F_3 = R$ induces maps $s_i: F_i @>>> F_{3-i}\spcheck$ and
a commutative diagram:
\begin{equation}
\label{dual}
\begin{CD}
\bold{P}^*: & \qquad & 0 @>>> R @>{d_3}>> F_2 @>{d_2}>> F_1
@>{d_1}>> R \\
&&&& @|  @VV{s_2}V  @VV{s_1}V  @| \\
(\bold{P}^*)\spcheck : && 0 @>>> R @>{d_1\spcheck}>> F_1\spcheck
@>{-d_2\spcheck}>> F_2\spcheck @>{d_3\spcheck}>> R
\end{CD}
\end{equation}
This map of complexes is an extension of the Gorenstein duality
isomorphism $R/I \cong \omega _{R/I} = \Ext ^3_R(R/I,R)$ to the
minimal projective resolutions of $R/I$ and $\omega _{R/I}$.  Since
any map between minimal projective resolutions which extends an
isomorphism in degree $0$ must be an isomorphism, it follows that the
$s_i$ are all isomorphisms.

We can therefore use the identification $s_2 : F_2\cong F_1\spcheck$.
A very simple computation (p.~465) shows that with this
identification, the commutativity and associativity of the
differential graded algebra structure on $\bold{P}^*$ imply the
skew-symmetry of $d_2$.  In particular $d_2$ must have even rank (say
$2p$), and $F_2$ must have odd rank $2p+1$.  The identification of
$d_1$ and $d_3$ with the vectors of Pfaffians of order $2p$ of $d_2$
is a lengthy but unproblematic computation (pp.~458--464).
\end{pf*}

Now let $X$ be a locally Gorenstein subcanonical subscheme of
codimension $3$ in $\Pnt$ with $\omega_X\cong\cal O_X (l)$.  We wish
to repeat the proof we have just sketched only with $\bold{P}^*$
replaced by a locally free resolution of $\cal O_X$:
\begin{equation}
\label{P}
\cal{P}^*: \qquad\quad 0 @>>> \cal{L} @>{d_3}>> \cal{F}_2
@>{d_2}>> \cal{F}_1 @>{d_1}>> \cal O_{\Pnt}
\end{equation}
where we will write $\cal{L}$ in place of $\omega_{\Pnt}(-l)$ in order
to simplify our diagrams.

A careful reading yields only two places where the fact that
$\bold{P}^*$ is a minimal projective resolution of $R/I$ was used in a
way that does not immediately carry over to the locally free
resolution $\cal{P}^*$.  The first place was in the definition of the
map $\Phi : S_2(\bold{P}^*) @>>> \bold{P}^*$ which made $\bold{P}^*$
into a commutative differential graded algebra.  Therefore we will
need to show directly the existence of a map of complexes
\[
\begin{CD}
S_2(\cal{P}^*): & \quad & \dotsb & \:\longrightarrow\: & \cal{L}
\oplus [\cal{F}_2\otimes \cal{F}_1] &
\:\stackrel{\sigma}{\longrightarrow}\: &
\cal{F}_2 \oplus \Lambda^2 \cal{F}_1 & \:\longrightarrow\: &
\cal{F}_1 & \:\longrightarrow\: & \cal O_{\Pnt} \\
&&&&@VV{\phi _3}V @VV{\phi _2}V @| @| \\
\cal{P}^*: && 0 & \:\longrightarrow\: & \cal{L} &
\:\stackrel{d_3}{\longrightarrow}\: & \cal{F}_2 &
\:\stackrel{d_2}{\longrightarrow}\: & \cal{F}_1 &
\:\stackrel{d_1}{\longrightarrow}\: & \cal O_{\Pnt}
\end{CD}
\]
The critical problem in defining the morphism of complexes is the
following.  Let $\psi : \Lambda^2 \cal{F}_1 @>>> \ker (d_1)$ be
defined by $\psi (a\wedge b)=d_1(a)b-d_1(b)a$. We then must lift
\begin{equation}
\label{liftdiag}
\begin{CD}
&&&&&&  \Lambda^2 \cal{F}_1 \\
&&&&&& @VV{\psi}V \\
0 @>>> \cal{L} @>>> \cal{F}_2 @>>> \ker (d_1) @>>> 0
\end{CD}
\end{equation}
to a $\phi\in\Hom (\Lambda ^2\cal{F}_1,\cal{F}_2)$. Once that is
done, the rest of the chain map follows.  For one may define $\phi _2
= (1_{\cal{F}_2},\phi )$.  Then
\[
\phi _2\circ\sigma (\cal{L}\oplus [\cal{F}_2\otimes \cal{F}_1])
\subset \ker (d_2) = \cal{L}.
\]
So $\phi _2\circ\sigma$ factors through $\cal{L}$, allowing one to
define $\phi _3$.  Thus one can put a commutative associative
differential graded algebra structure on $\cal{P}^*$ provided $\psi$
can be lifted. The obstruction to lifting $\psi$ lies in $\Ext
^1(\Lambda^2 \cal{F}_1,\cal{L}) \cong H^{n+2}(\Lambda^2
\cal{F}_1(l))^*$.

Once we have the commutative differential graded algebra structure on
$\cal{P}^*$, we may use it to define maps $s_i : \cal{F}_i @>>>
\cal{F}_{3-i}\spcheck \otimes \cal{L}$ and a commutative diagram
analogous to (\ref{dual}):
\begin{equation}
\label{dual2}
\begin{CD}
\cal{P}^*: & \qquad & 0 @>>> \cal{L} @>{d_3}>> \cal{F}_2 @>{d_2}>>
\cal{F}_1 @>{d_1}>> \cal O_{\Pnt} \\
&&&& @|  @VV{s_2}V  @VV{s_1}V  @| \\
(\cal{P}^*)\spcheck : &&0 @>>> \cal{L} @>{d_1\spcheck}>>
\cal{F}_1\spcheck\otimes\cal{L} @>{-d_2\spcheck}>>
\cal{F}_2\spcheck\otimes\cal{L} @>{d_3\spcheck}>> \cal O_{\Pnt}
\end{CD}
\end{equation}
The vertical maps extend the isomorphism $\cal O_X\cong \omega_X(-l)=
\cal{E}xt^3(\cal O_X,\cal{L})$.  We now run into the second problem
with locally free resolutions.  Namely, a morphism of locally free
resolutions which extends an isomorphism in degree $0$ is not
automatically an isomorphism between the resolutions.  But we reach
the conclusion:

\begin{proposition}
\label{conditions}
Suppose $X$ is a locally Gorenstein subcanonical subscheme of
codimension $3$ in $\Pnt$ with $\omega_X\cong\cal O_X (l)$. Then $X$
will be a Pfaffian scheme if $\cal O_X$ has a locally free resolution
$\cal{P}^*$ as in \rom{(}\ref{P}\rom{)} satisfying the following two
conditions:

\rom(a\rom) Any morphism of complexes $\cal{P}^* @>>>
(\cal{P}^*)\spcheck$ as in \rom{(}\ref{dual2}\rom{)} which extends the
identity of $\cal O_X$ is an isomorphism of complexes, and

\rom(b\rom) The morphism $\psi$ of \rom(\ref{liftdiag}\rom) lifts to a
map $\phi\in\Hom (\Lambda ^2\cal{F}_1,\cal{F}_2)$.
\end{proposition}

We will now construct locally free resolutions $\cal{P}^*$ satisfying
the conditions of the proposition.  Our method involves the Horrocks
correspondence.

\section{The Horrocks Correspondence}

In this section we give a modern description of the Horrocks
correspondence of \cite{H} using derived categories.  We include a
full proof of the principal properties of the correspondence from this
point of view (Theorem \ref{Horrocks}).  Taking advantage of the
greater flexibility of the derived category viewpoint, we develop a
technique which allows us to transfer a prescribed portion of the
cohomology of $\cal O_X$ to prescribed parts of a locally free
resolution (Proposition \ref{functorial}).

\subsection{Notation and Generalities}

We first recall some generalities about complexes.  If $\frak{A}$ is
an abelian category, let $C(\frak{A})$ (resp.\ $K(\frak{A})$,
$D(\frak{A})$) denote the category (resp.\ homotopy category, derived
category) of complexes of objects of $\frak{A}$, and let
$C^b(\frak{A})$, $C^-(\frak{A})$, $C^+(\frak{A})$, etc., denote the
corresponding complexes of bounded (resp.\ bounded above, bounded
below) complexes of objects of $\frak{A}$.  When speaking of
complexes, we will generally reserve the word ``isomorphism'' for
isomorphisms in $C(\frak A)$.  Isomorphisms in $K(\frak A)$ (resp.\
$D(\frak A)$) are referred to as homotopy equivalences (resp.\
quasi-isomorphisms).

If $r$ is an integer, then any complex $C^*$ of objects of $\frak{A}$
has two {\em canonical truncations} at $r$ and a {\em naive
truncation}:
\begin{align*}
\begin{CD}
\tau _{\leq r}(C^*): & \qquad & \cdots & \:\rightarrow\: & C^{r-2} &
\:\rightarrow\: & C^{r-1} & \:\rightarrow\: &
\ker(\delta^r) & \:\rightarrow\: & 0 & \:\rightarrow\: & 0 &
\:\rightarrow\: & \cdots , \\
\tau _{> r}(C^*): & \qquad & \cdots & \:\rightarrow\: & 0 &
\:\rightarrow\: & 0 & \:\rightarrow\: & C^r/\ker(\delta^r) &
\:\rightarrow\: & C^{r+1} & \:\rightarrow\: & C^{r+2} &
\:\rightarrow\: & \cdots . \\
\sigma _{\geq r}(C^*): & \qquad & \cdots & \:\rightarrow\: & 0 &
\:\rightarrow\: & 0 & \:\rightarrow\: & C^r & \:\rightarrow\: &
C^{r+1} & \:\rightarrow\: & C^{r+2} & \:\rightarrow\: & \cdots .
\end{CD}
\end{align*}
All the truncations are functorial in $C(\frak{A})$. The canonical
truncations are functorial in $K(\frak{A})$ and $D(\frak{A})$ as well.
We will often find it more convenient to write $\tau_{<r+1}$ instead
of $\tau_{\leq r}$.

Suppose now that $\frak{A}$ has enough projectives.  Every bounded
above complex $C^*$ of objects in $\frak{A}$ admits a {\em projective
resolution}, i.e.\ a quasi-isomorphism $P^* @>>> C^*$ with $P^*$ a
complex of projectives (\cite{Ha} Proposition I.4.6).  The projective
resolution of a complex is unique up to homotopy equivalence.  If
$C^*$ and $E^*$ are bounded above complexes of objects in $\frak{A}$,
and if $P^* @>>> C^*$ is a projective resolution of $C^*$, then there
is a natural isomorphism $\Hom_{D^-(\frak{A})}(C^*,E^*) \cong
\Hom_{K^-(\frak{A})}(P^*,E^*)$.  In particular if $\frak{P}$ denotes
the full subcategory of projective objects of $\frak{A}$, then the
natural functor $K^-(\frak{P}) @>>> D^-(\frak{A})$ is an equivalence
of categories (\cite{Ha} Proposition I.4.7).  This can be refined to
the following statement:

\begin{lemma}
\label{category}
Let $\frak{A}$ be an abelian category with enough projectives, and let
$\frak{P}$ be the full subcategory of projective objects of
$\frak{A}$.  Suppose $A\subset D^-(\frak{A})$ and $P\subset
K^-(\frak{P})$ are full subcategories such that $\ob(P)\subset \ob(A)$
and every object of $A$ has a projective resolution belonging to $P$.
Then the natural functor $P @>>> A$ is an equivalence of categories.
\end{lemma}

Let $S=k[X_0,\dots,X_N]$ be the homogeneous coordinate ring of $\PN$,
and let $\frak{m}=(X_0,\dots,X_N)$ be its irrelevant ideal.  Let
$\grsmod$ be the category of graded $S$-modules.  Then $\grsmod$ has
enough projectives, namely the free modules.  We will call a complex
$P^*$ of projectives in $\grsmod$ a {\em minimal} if all its objects
$P^i$ are free of finite rank and its differential $\delta^*$
satisfies $\delta^i(P^i) \subset \frak m P^{i+1}$ for all $i$.  If
$C^*$ is a bounded above complex of objects in $\grsmod$ whose
cohomology modules $H^i(C^*)$ are all finitely generated, then $C^*$
has a {\em minimal projective resolution}, i.e.\ a projective
resolution by a minimal complex of projectives.  The next lemma, which
is a well known consequence of Nakayama's lemma, says that minimal
projective resolutions are unique up to isomorphism and not merely up
to homotopy equivalence:

\begin{lemma}
\label{NAK}
Let $\phi: P^* @>>> Q^*$ be a homotopy equivalence between {\em
minimal} complexes of free graded $S$-modules of finite rank.  Then
$\phi$ is an isomorphism.
\end{lemma}


Let $\OMod$ be the category of sheaves of $\cal{O}_{\PN}$-modules.
For $\cal{E}$ a sheaf of $\cal{O}_{\PN}$-modules, let $\Gamma
_*(\cal{E})=\bigoplus _{t\in\Bbb{Z}} \Gamma(\cal{E}(t))$.  Then
$\Gamma_*$ defines a left exact functor from $\OMod$ to $\grsmod$.  It
has a right derived functor $\bold{R}\Gamma_*: D^b(\OMod) @>>>
D^b(\grsmod)$ whose cohomology functors we denote $H^i_*(\cal{E}) =
\bigoplus _{t\in\Bbb{Z}}H^i(\cal{E}(t))$.  The functor $\Gamma_*$ has
an exact left adjoint $\widetilde{\ \ }$, the functor of associated
sheaves.

Let $\Gamma_{\frak m}: \grsmod @>>> \grsmod$ be the functor
associating to a graded $S$-module $M$ the maximal submodule
$\Gamma_{\frak m}(M) \subset M$ supported at the origin $0$ of $\Bbb
A^{N+1}$.  This functor is also left exact and has a right derived
functor $\bold R\Gamma_{\frak m}: D^b(\grsmod) @>>> D^b(\grsmod)$.
Its cohomology functors are denoted $H^i_{\frak m}$.

\begin{lemma}
\label{bounds}
Let $P^*$ be a bounded complex of free graded $S$-modules of finite
rank where $S=k[X_0,\dots,X_N]$.  If $P^*$ is minimal, then
\begin{align*}
\max\{i\mid P^i\neq 0\} = & \max\{i\mid H^i(P^*)\neq 0\},\\
\min\{i\mid P^i\neq 0\} = & \min\{i\mid H^i_{\frak m}(P^*)\neq 0\}
-N-1.
\end{align*}
\end{lemma}

\begin{pf}
The assertion about maxima is a simple and well-known application of
the minimality condition and Nakayama's Lemma.  The assertion about
minima, which is essentially the Auslander-Buchsbaum theorem, reduces
to the assertion about maxima by Serre duality.
\end{pf}

\subsection{The Horrocks Correspondence}

We now begin to describe the components of the Horrocks
correspondence.  Let $\frak{B}$ be the full subcategory of $\OMod$ of
locally free sheaves of finite rank, and let $\frak{Z}$ denote the
full category of $D^b(\grsmod)$ of complexes $C^*$ such that
$H^i(C^*)$ is of finite length for $0<i<N$ and $H^i(C^*)$ vanishes for
all other $i$.

The Horrocks correspondence consists of a functor $\zeta: \frak{B}
@>>> \frak{Z}$ and a map $\cal{H}: \ob(\frak{Z}) @>>> \ob(\frak{B})$
in the opposite direction.  The functor $\zeta$ is simply
$\tau_{>0}\tau_{< N}\bold{R}\Gamma _*$.  For $\cal{E}$ a vector bundle
on $\PN$, the cohomology of $\zeta(\cal{E})$ is of course:
\[
H^i(\zeta(\cal{E}))=
\begin{cases}
H^i_*(\cal{E}) & \text{if }0<i<N,\\
0 & \text{otherwise.}
\end{cases}
\]
Since $\cal{E}$ is locally free of finite rank, $H^i_*(\cal{E})$ is of
finite length for $0<i<N$.  So $\zeta(\cal{E})\in\ob(\frak{Z})$.

We now define $\cal{H}$.  Any $C^*\in\ob(\frak{Z})$ has a minimal
projective resolution $P^* @>>> C^*$.  We define $\cal{H}(C^*)$ to be
the kernel of the differential $\widetilde\delta^0:\widetilde P^0 @>>>
\widetilde P^1$.  Then $\cal{H}(C^*)$ is a vector bundle because it
fits into an exact complex of vector bundles
\begin{equation}
\label{Hcomp}
\dotsb @>>> 0 @>>> \cal{H}(C^*) @>>> \widetilde P^0 @>>> \widetilde
P^1 @>>> \dotsb @>>> \widetilde P^{N-1} @>>> 0 @>>> \dotsb .
\end{equation}
Note that $\cal{H}(C^*)$ is well-defined up to isomorphism because the
minimal projective resolution $P^*$ of $C^*$ is unique up to
isomorphism because of Lemma \ref{NAK}.  However, $\cal{H}$ is not a
functor.

The principal results of Horrocks' paper \cite{H} can be described in
the following way:

\begin{theorem}[Horrocks]
\label{Horrocks}
Let $\frak{B}$ be the category of locally free sheaves of finite rank
on $\PN$, and let $\frak{Z}$ be the full subcategory of $D^b(\grsmod)$
of complexes $C^*$ such that $H^i(C^*)$ is of finite length if
$0<i<N$, and $H^i(C^*)=0$ for all other $i$.  Let $\zeta =
\tau_{>0}\tau_{< N}\bold{R}\Gamma _*: \frak{B} @>>> \frak{Z}$, and let
$\cal{H}:\ob(\frak{Z}) @>>> \ob(\frak{B})$ be the map defined as in
\rom(\ref{Hcomp}\rom) above.

\rom(a\rom) If $\cal{E}\in\ob(\frak{B})$, then $\cal{E}
\cong \cal{H}\zeta(\cal{E}) \oplus \bigoplus _i
\cal{O}_{\PN}(n_i)$ for some integers $n_i$.

\rom(b\rom) If $C^*\in\ob(\frak{Z})$, then  $\zeta\cal{H}(C^*)\simeq
C^*$.

\rom(c\rom) If $\cal{E},\cal{F}\in\ob(\frak{B})$, then
$\Hom_{\frak{Z}}(\zeta(\cal{E}),\zeta(\cal{F}))\cong
\Hom(\cal{E},\cal{F})/\Hom_{\Phi}(\cal{E},\cal{F})$ where
$\Hom_{\Phi}(\cal{E},\cal{F})$ is the set of all morphisms which
factor through a direct sum of line bundles.
\end{theorem}

The theorem may be read as saying the following. Call two vector
bundles $\cal{E}$ and $\cal{F}$ {\em stably equivalent} if there exist
sets of integers $\{n_i\}$ and $\{m_j\}$ such that $\cal{E} \oplus
\bigoplus_i\cal{O}_{\PN}(n_i) \cong \cal{F}\oplus \bigoplus_j
\cal{O}_{\PN}(m_j)$.  Then the theorem says that $\zeta$ and $\cal{H}$
induce a one-to-one correspondence between stable equivalence classes
of vector bundles on $\PN$ and quasi-isomorphism classes of complexes
in $\frak{Z}$.

For Horrocks' proof of the theorem, see \cite{H} Lemma 7.1 and Theorem
7.2 and the discussion between them.  However, Horrocks' definition of
the category $\frak{Z}$ and the functor $\zeta$ are different from
ours, and demonstrating the equivalence of the definitions is somewhat
tedious.  So instead of referring the reader to Horrocks' paper, we
give a new proof.  The first step is the following lemma:

\begin{lemma}
\label{tauresol}
\rom(a\rom) Suppose
\[
P^*:\qquad \dotsb @>>> 0 @>>> P^0 @>>> P^1 @>>> \dotsb @>>>
P^{N-1} @>>> 0 @>>> \dotsb
\]
is a complex of free graded $S$-modules of finite rank such that
$H^i(P^*)$ is a module of finite length for $0<i<N$.  Let $\cal{E}=
H^0(P^*)\sptilde$.  Then $P^*$ is quasi-isomorphic to
$\tau_{<N}\bold{R}\Gamma_*(\cal{E})$.

\rom(b\rom) Conversely, if $\cal{E}$ is a vector bundle on $\PN$, then
the minimal projective resolution of
$\tau_{<N}\bold{R}\Gamma_*(\cal{E})$ is of the above form.
\end{lemma}

\begin{pf}
(a) Note that the complex $\widetilde{P}^*$ of coherent sheaves on
$\PN$ has vanishing cohomology in degrees different from $0$.  So it
is quasi-isomorphic to $H^0(\widetilde{P}^*) = \cal{E}$.  Hence the
triangle of functors of \cite{W} Proposition 1.1:
\[
\bold{R}\Gamma_{\frak m} @>>> \Id @>>> \bold{R}\Gamma_*\circ\sptilde
@>>> \bold{R}\Gamma_{\frak m}[1],
\]
when applied to $P^*$, yields a triangle
\begin{equation}
\label{triangle}
\bold{R}\Gamma_{\frak m}(P^*) @>>> P^* @>\beta>>
\bold{R}\Gamma_*(\cal{E}) @>>> \bold{R}\Gamma_{\frak m}(P^*)[1].
\end{equation}
By Lemma \ref{bounds}, we have $H^i_{\frak m}(P^*)=0$ for $i \leq N$.
So $H^i(\beta): H^i(P^*) @>>> H^i_*(\cal{E})$ is an isomorphism for
$i<N$.  Therefore $\beta$ induces a quasi-isomorphism of $P^*$ onto
$\tau_{<N}\bold{R}\Gamma_* (\cal{E})$.

(b) Conversely, if $\cal{E}$ is a vector bundle on $\PN$, then
$H^i_*(\cal{E})$ is finitely generated for $i<N$.  Hence
$\tau_{<N}\bold{R}\Gamma_*(\cal{E})$ has a minimal projective
resolution $P^*$.  For $0<i<N$ the module $H^i(P^*) = H^i_*(\cal E)$
is of finite length because $\cal E$ is locally free.  By construction
$H^i(P^*) = H^i(\tau_{<N}\bold R\Gamma_*(\cal E)) = 0$ for $i \geq N$.
So we have $P^i = 0$ for $i\geq N$ by Lemma \ref{bounds}.  Looking
again at the triangle (\ref{triangle}), we see by the construction of
$P^*$ that $H^i(\beta)$ is an isomorphism for $i<N$ and an injection
for $i=N$.  So $H^i_{\frak m}(P^*) = 0$ for $i\leq N$.  So by Lemma
\ref{bounds} we see that $P^i=0$ for $i\leq -1$.  Thus $P^*$ has the
form asserted by the lemma.
\end{pf}

We now wish to functorialize the previous lemma.  Let $B\subset
K^b(\grsmod)$ be the full subcategory of complexes of the form
\begin{equation}
\label{rightresol}
\dotsb @>>> 0 @>>> P^0 @>>> P^1 @>>> \dotsb @>>> P^{N-1} @>>>
0 @>>> \dotsb
\end{equation}
such that the $P^i$ are free of finite rank for all $i$, the modules
$H^i(P^*)$ are of finite length for $0<i<N$ and the differentials
satisfy $\delta^i(P^i)\subset \frak{m}P^{i+1}$ for all $i$.  For any
vector bundle $\cal E$ on $\PN$ we now define $P^*(\cal E)$ as the
minimal projective resolution of $\tau_{<N}\bold R\Gamma_* (\cal E)$.
By Lemma \ref{tauresol}, $P^*(\cal E)$ is always an object of $B$.

\begin{lemma}
\label{Bequiv}
The functor $P^*: \frak{B} @>>> B$ which associates to an $\cal{E}
\in\ob(\frak{B})$ the minimal projective resolution of $\tau_{<N}
\bold{R}\Gamma_*(\cal{E})$ is an equivalence of categories with
inverse given by $C^*\mapsto H^0(C^*)\sptilde$.
\end{lemma}

\begin{pf}
Since the functor $\tau_{<N}\bold R\Gamma_*: \frak B @>>>
D^-(\grsmod)$ has a left inverse $H^0(-)\sptilde$, it induces an
equivalence between $\frak B$ and the full subcategory $A\subset
D^-(\grsmod)$ of complexes quasi-isomorphic to complexes in the image
of $\tau_{<N}\bold R\Gamma_*$.  But by Lemma \ref{tauresol}, the full
subcategory $B\subset K^-(\grsmod)$ has the properties that
$\ob(B)\subset\ob(A)$ and that the minimal projective resolution of
every object of $A$ belongs to $B$.  Hence the natural functor $B @>>>
A$ is also an equivalence of categories by Lemma \ref{category}.
Since $P^*$ is exactly the composition of the equivalence
$\tau_{<N}\bold R\Gamma_*: \frak B @>>> A$ with the inverse of the
equivalence $B @>>> A$, it is an equivalence.  The inverse of $P^*$
remains the same as that of $\tau_{<N} \bold R\Gamma_*$, namely
$H^0(-)\sptilde$.
\end{pf}

Now the graded module associated to a vector bundle $\cal E$ on $\PN$
has a minimal projective resolution:
\[
0 @>>> Q^{-(N-1)} @>>> \dotsb @>>> Q^{-1} @>>> Q^0 @>>> \Gamma_*(\cal
E)
\]
For any $\cal E$ we now define the following complexes in addition to
the $P^*(\cal E)$ defined above.  First we set:
\[
Q^*(\cal E):\qquad \dotsb @>>> 0 @>>> Q^{-(N-1)} @>>> \dotsb @>>>
Q^{-1} @>>> Q^0 @>>> 0 @>>> \dotsb.
\]
We then let $R^*(\cal E)$ be the natural concatenation of $Q^*(\cal
E)$ with $P^*(\cal E)$ induced by the composition $Q^0
\twoheadrightarrow \Gamma_*(\cal E) \hookrightarrow P^0$:
\[
R^*(\cal E):\qquad\dotsb @>>> 0 @>>> Q^{-(N-1)} @>>> \dotsb @>>> Q^0
@>>> P^0 @>>> \dotsb @>>> P^{N-1} @>>> 0 @>>> \dotsb
\]
Thus $R^i(\cal E) = P^i(\cal E)$ for $i\geq 0$, and $R^i(\cal E) =
Q^{i+1}(\cal E)$ for $i<0$.  Note that although the projective
complexes $P^*(\cal E)$ and $Q^*(\cal E)$ are minimal, $R^*(\cal E)$
may not be minimal, because there may be a direct factor of $Q^0(\cal
E)$ which is mapped isomorphically onto a direct factor of $P^0(\cal
E)$.  However, one may write $R^*(\cal E)$ as the direct sum of a
minimal complex of projectives $R^*_{\min}(\cal E)$
\begin{align*}
R^*_{\min}(\cal E): \qquad \dotsb @>>> Q^{-2} @>>> Q^{-1}  @>>>
Q^0_{\min} & @>>> P^0_{\min} @>>> P^1 @>>> P^2 @>>> \dotsb\\
\intertext{and of an exact complex of projectives}
\stepcounter{equation}\tag{\theequation}\label{L:module}
\dotsb @>>> 0 @>>> L & @>\Id>> L @>>> 0 @>>> \dotsb.
\end{align*}

The complexes $Q^*(\cal E)$, $R^*(\cal E)$, and $R^*_{\min}(\cal E)$
are all functorial (in the homotopy category) in $\cal E$.  Moreover,
we may use the identification between the categories $\frak B$ and $B$
to define complexes $Q^*(P^*)$, $R^*(P^*)$, and $R^*_{\min}(P^*)$ for
$P^*$ in $B$.  Namely, $Q^*(P^*)$ is the minimal projective resolution
of $H^0(P^*)$, $R^*(P^*)$ is the concatenation of $Q^*(P^*)$ with
$P^*$, etc.

We now define a homotopy category of complexes of type $R^*_{\min}$.
More formally, let $Z\subset K^b(\grsmod)$ be the full subcategory of
minimal complexes of projective modules of finite rank of the form
\begin{equation}
\label{R:complex}
\dotsb @>>> 0 @>>> R^{-N} @>>> \dotsb @>>> R^{-1} @>>> R^0 @>>> \dotsb
@>>> R^{N-1} @>>> 0 @>>> \dotsb
\end{equation}
such that the cohomology modules $H^i(R^*)$ are of finite length for
$0<i<N$ and vanish for all other $i$.

We need one more lemma before proving Theorem \ref{Horrocks}.

\begin{lemma}
\label{Zequiv}
The natural functor $Z @>>> \frak{Z}$ is an equivalence of categories.
\end{lemma}

\begin{pf}
Let $R^*$ be the minimal projective resolution of an object $C^*$ of
$\frak Z$.  Since $H^i(R^*)= H^i(C^*)=0$ for $i\geq N$, we have
$R^i=0$ for $i\geq N$ by Lemma \ref{bounds}.  Moreover, all the
$H^i(C^*)$ are of finite length, so $H^i_{\frak m}(C^*) = H^i(C^*)$
for all $i$.  In particular, $H^i_{\frak m}(R^*)=H^i_{\frak m}(C^*)
=0$ for $i\leq 0$.  So $R^i=0$ for $i\leq -N-1$ by Lemma \ref{bounds}.
Thus the minimal projective resolution of any object of $\frak Z$ is
in $Z$.  The lemma now follows from Lemma \ref{category}.
\end{pf}

\begin{pf*}{Proof of Theorem \ref{Horrocks}}
Lemma \ref{Bequiv} permits us to identify a vector bundle $\cal E$
with the complex $P^*(\cal E)$ of $B$.  Since $P^*(\cal E)$ is already
quasi-isomorphic to $\tau_{<N}\bold R\Gamma_*(\cal E)$, the complex
$\zeta(\cal E) = \tau_{>0}\tau_{<N} \bold R\Gamma_*(\cal E)$ is
quasi-isomorphic to the complex
\[
\dotsb @>>> 0 @>>> \Gamma_*(\cal E) @>>> P^0 @>>> P^1 @>>> \dotsb @>>>
P^{N-1} @>>> 0 @>>> \dotsb
\]
and hence to the complexes $R^*(\cal E)$ and $R^*_{\min}(\cal E)$.
Hence the object $\zeta(\cal E)$ in $\frak Z$ is quasi-isomorphic to
the object $R^*_{\min}(\cal E)$ of $Z$.  Hence after identifying
$\frak B$ and $\frak Z$ with $B$ and $Z$ by Lemmas \ref{Bequiv} and
\ref{Zequiv}, the functor $\zeta$ may be identified with the functor
from $B$ to $Z$ which associates to any complex $P^*$ in $B$ the
corresponding complex $R^*_{\min}$ as described earlier.

Similarly, given any object $C^*$ of $\frak Z$ with minimal projective
resolution $R^*$, the definitions say that $P^*(\cal H(C^*))=
\sigma_{\geq 0}(R^*)$, the naive truncation.  Thus the map $\cal H:
\ob(\frak Z) @>>> \ob(\frak B)$ may be identified with $\sigma_{\geq
0}: \ob(Z) @>>> \ob(B)$.  Note that since all objects of $Z$ and $B$
are minimal complexes of projective modules, homotopy equivalence
classes of objects of $Z$ and $B$ coincide with isomorphism classes.
Hence the map $\sigma_{\geq 0}: \ob(Z) @>>> \ob(B)$ preserves homotopy
equivalence.  Since $Z$ and $B$ are subcategories of the homotopy
category, this means that $\sigma_{\geq 0}$ is well-defined on objects
of $Z$.  However, $\sigma_{\geq 0}$ and hence $\cal H$ are not
well-defined on morphisms of $Z$.

(a) The above identifications now say if $\cal E\in\ob(\frak B)$, then
$\cal H\zeta(\cal E)$ is the object of $\frak B$ corresponding to the
complex $\sigma_{\geq 0}(R^*_{\min}(\cal E))$:
\[
\sigma_{\geq 0}(R^*_{\min}(\cal E)): \qquad \dotsb @>>> 0 @>>>
P^0_{\min} @>\mu>> P^1 @>>> \dotsb @>>> P^{N-1} @>>> 0 @>>> \dotsb .
\]
By Lemma \ref{Bequiv}, the sheaf $\cal H\zeta(\cal E)$ is
$\ker(\mu)\sptilde$.  So $\cal E = \cal H\zeta(\cal
E)\oplus\widetilde{L}$ where $L$ is the projective module of
(\ref{L:module}).  Since $\widetilde{L}$ is now a direct sum of line
bundles, (a) follows.

(b) If $C^*$ is an object of $\frak Z$ with minimal projective
resolution $R^*$ in $Z$ of the form (\ref{R:complex}), then the above
computations identify $\cal H(C^*)$ in $\frak B$ with $P^*(\cal
H(C^*)) = \sigma_{\geq 0}(R^*)$ in $B$.  Thus $\zeta\cal H(C^*)$
becomes identified with $R^*_{\min}(\cal H(C^*))$ which is just $R^*$
again.  Since $R^*$ is quasi-isomorphic to $C^*$, we have $\zeta\cal
H(C^*)\simeq C^*$ as desired.

(c) After identifying $\frak B$ with $B$ and $\frak Z$ with $Z$,
assertion (c) becomes the statement: For any pair of objects $E^*$ and
$F^*$ in $B$, the natural map
\begin{equation}
\label{B:to:Z}
\Hom_B(E^*,F^*) @>>> \Hom_Z(R^*_{\min}(E^*),R^*_{\min}(F^*))
\end{equation}
is surjective and its kernel is the subspace of morphisms which factor
through an object of $B$ of the form
\begin{equation}
\label{triv:comp}
\dotsb @>>> 0 @>>> L @>>> 0 @>>> \dotsb
\end{equation}
with $L$ a free graded $S$-module of finite rank appearing in degree
$0$.

We first prove surjectivity.  Suppose $\phi\in \Hom_Z
(R^*_{\min}(E^*), R^*_{\min}(F^*))$.  Since $Z$ is a homotopy
category, $\phi$ is actually a homotopy equivalence class of maps in
$C(\grsmod)$.  So we may choose a chain map $f$ in the class $\phi$.
Then $f$ may be extended to a chain map $\overline f: R^*(E^*) @>>>
R^*(F^*)$ by defining it to be $0$ on the exact factor of the type
(\ref{L:module}). Then $\sigma_{\geq 0}\overline f$ maps $E^*$ to
$F^*$, and its homotopy class in $B$ has image $\phi$ in $Z$.  This
proves surjectivity.

We now compute the kernel of (\ref{B:to:Z}).  First if $\alpha\in
\Hom_B(E^*,F^*)$ factors through a complex $L^*$ of the form
(\ref{triv:comp}), then $R^*_{\min}(\alpha)$ factors through
$R^*_{\min}(L^*)=0$ and so vanishes.  So the kernel of (\ref{B:to:Z})
contains all morphisms which factor through complexes of the form
(\ref{triv:comp}).

Conversely, suppose $\alpha$ is in the kernel of (\ref{B:to:Z}).
Since $\alpha$ is a morphism in $B$, it is a homotopy class of chain
maps from which we may choose a member $\beta$.  We may complete
$\beta$ to a chain map $\rho: R^*(E^*) @>>> R^*(F^*)$.
\[
\begin{CD}
R^*(E^*) & \qquad & \dotsb & \:\rightarrow\: & 0 & \:\rightarrow\: &
\overline E^{-N} & \:\rightarrow\: & \dotsb & \:\rightarrow\: &
\overline E^{-1} & \:\rightarrow\: & E^0 & \:\rightarrow\: & \dotsb &
\:\rightarrow\: & E^{N-1} & \:\rightarrow\: & 0 & \:\rightarrow\: &
\dotsb \\
@VV{\rho}V && @VVV @VVV && @VVV @VV{\beta}V && @VV{\beta}V @VVV \\
R^*(F^*) & \qquad & \dotsb & \:\rightarrow\: & 0 & \:\rightarrow\: &
\overline F^{-N} & \:\rightarrow\: & \dotsb & \:\rightarrow\: &
\overline F^{-1} & \:\rightarrow\: & F^0 & \:\rightarrow\: & \dotsb &
\:\rightarrow\: & F^{N-1} & \:\rightarrow\: & 0 & \:\rightarrow\: &
\dotsb
\end{CD}
\]
The homotopy class of $\rho$ is the image of $\alpha$ under $R^*$ and
so must vanish by hypothesis.  (Note that $R^*$ and $R^*_{\min}$ are
homotopy equivalent.)  Thus $\rho$ is homotopic to $0$.  Thus if we
write $\delta^i$ for the differentials of $R^*(E^*)$, and $\epsilon^i$
for the differentials of $R^*(F^*)$, then there is a chain homotopy $h
= (h^i)$ such that $\rho^i = h^{i+1}\delta^i + \epsilon^{i-1}h^i$ for
all $i$.  Now restrict $h$ to a chain homotopy $\widehat{h} =
(\widehat{h}^i)$ with $\widehat{h}^i: E^i @>>> F^{i-1}$ defined by
defined by $\widehat{h}^i = h^i$ for all $i\geq 1$, and $\widehat{h}^i
= 0$ for all $i\leq 0$.  Then $\beta$ is homotopic to a morphism whose
components are
\[
\beta^i-(\widehat{h}^{i+1}\delta^i + \epsilon^{i-1}\widehat{h}^i)=
\begin{cases}
\rho^i-(h^{i+1}\delta^i+\epsilon^{i-1}h^i) = 0 & \text{if }i\geq1,\\
\rho^0-h^1\delta^0 = \epsilon^{-1}h^0 & \text{if }i=0,\\
0 & \text{if }i\leq -1.
\end{cases}
\]
Hence the homotopy class $\alpha$ of $\beta$ factors through the complex
\[
\dotsb @>>> 0 @>>> \overline F^{-1} @>>> 0 @>>> \dotsb
\]
of type (\ref{triv:comp}).  So the kernel of (\ref{B:to:Z}) is as
asserted.  This completes the proof of the theorem.
\end{pf*}

We will use three further results concerning the Horrocks
correspondence.  The first will permit us to use the Horrocks
correspondence to constuct locally free resolutions of coherent
sheaves.

\begin{proposition}
\label{functorial}
Let $\cal{Q}$ be a quasi-coherent sheaf on $\PN$, let
$C^*\in\ob(\frak{Z})$, and let $\beta : C^* @>>> \tau_{>0}\tau_{<N}
\bold{R}\Gamma _*(\cal{Q})$ be a morphism in $D^b(\grsmod)$.  Then there
exists a morphism of quasi-coherent sheaves $\widetilde{\beta} :
\cal{H}(C^*) @>>> \cal{Q}$ such that $\beta =\tau_{>0}\tau_{<N}
\bold{R}\Gamma _*(\widetilde{\beta})$.  In particular, the induced
morphisms $H^i_*(\cal{H}(C^*)) @>>> H^i_*(\cal{Q})$ are the same as
$H^i(\beta)$ for $1\leq i\leq N-1$.
\end{proposition}

\begin{pf}
Let $R^*$ be a minimal projective resolution of $C^*$, and let
$\cal{I}^*$ be an injective resolution of $\cal{Q}$.  Then $\beta$ may
be identified with an actual chain map
\[
\begin{CD}
\dotsb & \:\rightarrow\: & R^{-2} & \:\rightarrow\: & R^{-1} &
\:\rightarrow\: & R^0 & \:\stackrel{\lambda}{\longrightarrow}\: & R^1 &
\:\rightarrow\: & \dotsb & \:\rightarrow\: & R^{N-2} & \:\rightarrow\:
& R^{N-1} & \:\rightarrow\: & 0 \\
&& @VVV @VVV @VVV @VVV && @VVV @VVV \\
\dotsb & \:\rightarrow\: & 0 & \:\rightarrow\: & \Gamma_*(\cal{Q}) &
\:\rightarrow\: & \Gamma_*(\cal{I}^0) & \:\stackrel{\mu}{\longrightarrow}\: &
\Gamma_*(\cal{I}^1) & \:\rightarrow\: & \dotsb & \:\rightarrow\: &
\Gamma_*(\cal{I}^{N-2}) & \:\rightarrow\: & \ker(\delta^{N-1}) &
\:\rightarrow\: & 0
\end{CD}
\]
Thus $\beta$ induces a morphism $\widetilde{\beta}$ from $\cal{H}(C^*)
= \ker(\lambda)\sptilde$ to $\cal{Q} = \ker(\mu)\sptilde$.

We now need to calculate $\bold{R}\Gamma_*(\widetilde{\beta})$.
Consider the complex
\[
P^*: \qquad \dotsb @>>> 0 @>>> R^0 @>>> R^1 @>>> \dotsb @>>> R^{N-1}
@>>> 0 @>>> \dotsb .
\]
The previous diagram induces a new commutative diagram
\[
\begin{CD}
\widetilde{P}^*&: & \qquad\qquad & \dotsb & \:\rightarrow\: & 0 &
\:\rightarrow\: & \widetilde R^0 & \:\rightarrow\: & \widetilde R^1 &
\:\rightarrow\: & \dotsb & \:\rightarrow\: & \widetilde R^{N-2} &
\:\rightarrow\: & \widetilde R^{N-1} & \:\rightarrow\: & 0 &
\:\rightarrow\: & 0 & \:\rightarrow\: & \dotsb \\
@VV{\overline \beta}V &&& @VVV @VVV @VVV && @VVV @VVV @VVV @VVV\\
\cal{I}^*&: && \dotsb & \:\rightarrow\: & 0 & \:\rightarrow\: &
\cal{I}^0 & \:\rightarrow\: & \cal{I}^1 & \:\rightarrow\: & \dotsb &
\:\rightarrow\: & \cal{I}^{N-2} & \:\rightarrow\: & \cal{I}^{N-1} &
\:\rightarrow\: & \cal{I}^N & \:\rightarrow\: & 0 & \:\rightarrow\: &
\dotsb
\end{CD}
\]
between resolutions of $\cal{H}(C^*)$ and $\cal{Q}$ extending
$\widetilde \beta$.  Let $\gamma: P^* @>>> J^*$ be an injective
resolution of $P^*$.  Then $\overline\beta$ factors through
$\widetilde\gamma$ as $\widetilde{P}^* @>>> \widetilde{J}^* @>>>
\cal{I}^*$. Applying $\Gamma_*$ now gives a factorization
\begin{equation}
\label{factorize}
P^* @>>> \Gamma_*(\widetilde{J}^*) @>>> \Gamma_*(\cal{I}^*).
\end{equation}
Now $\overline\beta$ is a map between resolutions of $\cal{H}(C^*)$
and $\cal{Q}$, respectively, which extends $\widetilde{\beta}:
\cal{H}(C^*) @>>> \cal{Q}$, while $\widetilde\gamma$ is a
quasi-isomorphism.  So the map $\widetilde{J}^* @>>> \cal{I}^*$ is a
map between injective resolutions of $\cal{H}(C^*)$ and $\cal{Q}$
extending $\widetilde\beta$.  So by definition, the second arrow of
(\ref{factorize}) is $\bold{R}\Gamma_*(\widetilde{\beta}) :
\bold{R}\Gamma_*(\cal{H}(C^*)) @>>> \bold{R}\Gamma_*(\cal{Q})$.  On
the other hand, the proof of Lemma \ref{tauresol}(a) shows that the
first arrow of (\ref{factorize}) can be identified with the truncation
$\tau_{<N}(\bold{R}\Gamma_*(\cal{H}(C^*))) @>>>
\bold{R}\Gamma_*(\cal{H}(C^*))$ because it induces isomorphisms
$H^i(P^*) \cong H^i(\Gamma_*(\widetilde{J}^*)) = H^i_*(\cal{H}(C^*))$
for $i < N$.  Hence $\Gamma_*(\overline\beta): P^* @>>>
\Gamma_*(\cal{I}^*)$ can be identified with the composition of the
truncation $\tau_{<N}(\bold{R}\Gamma_*(\cal{H}(C^*))) @>>>
\bold{R}\Gamma_*(\cal{H}(C^*))$ with $\bold{R}\Gamma_*
(\widetilde{\beta})$.  Thus $\tau_{<N}\bold{R}\Gamma_*
(\widetilde{\beta})$ may be identified with the diagram
\[
\begin{CD}
\dotsb & \:\rightarrow\: & 0 &
\:\rightarrow\: & R^0 & \:\rightarrow\: & R^1 &
\:\rightarrow\: & \dotsb & \:\rightarrow\: & R^{N-2} & \:\rightarrow\:
& R^{N-1} & \:\rightarrow\: & 0 & \:\rightarrow\: & \dotsb \\
&& @VVV @VVV @VVV && @VVV @VVV @VVV\\
\dotsb & \:\rightarrow\: & 0 &
\:\rightarrow\: & \Gamma_*(\cal{I}^0) & \:\rightarrow\: &
\Gamma_*(\cal{I}^1) & \:\rightarrow\: & \dotsb & \:\rightarrow\: &
\Gamma_*(\cal{I}^{N-2}) & \:\rightarrow\: & \ker(\delta^{N-1}) &
\:\rightarrow\: & 0 & \:\rightarrow\: & \dotsb
\end{CD}
\]
induced by $\beta$.  Truncating on the left, we reach a diagram
equivalent to the first diagram of the proof of the proposition.  So
$\beta = \tau_{>0}\tau_{<N} \bold{R}\Gamma _*(\widetilde{\beta})$.
\end{pf}

We now need two homological criteria for maps of vector bundles to be
isomorphisms.

\begin{lemma}
\label{isom}
Let $\cal E$ and $\cal F$ be vector bundles on $\PN$ with neither
containing a line bundle as a direct factor.  If $\alpha: \cal E @>>>
\cal F$ is a map such that $H^i_*(\alpha): H^i_*(\cal E) @>>>
H^i_*(\cal F)$ is an isomorphism for $0<i<N$, then $\alpha$ is an
isomorphism.
\end{lemma}

\begin{pf}
We use the notation of the proof of Theorem \ref{Horrocks}.  Let
$E^*=P^*(\cal E)$ and $F^*=P^*(\cal F)$, and let $\overline\alpha: E^*
@>>> F^*$ be the map induced by $\alpha$.  The hypothesis $\cal E =
\cal H\zeta(\cal E)$ implies that $E^*$ is homotopy equivalent to
$\sigma_{\geq 0}R^*_{\min}(E^*)$, or equivalently that $R^*(E^*)$ is a
minimal complex of projectives.  Similarly, $R^*(F^*)$ is a minimal
complex of projectives.  The hypothesis on $\alpha$ implies that
$\zeta(\alpha): \zeta(\cal E) @>>> \zeta(\cal F)$ is a
quasi-isomorphism.  This in turn translates into
$R^*_{\min}(\overline\alpha)$ being a homotopy equivalence.  But
because of the earlier hypotheses, this means that
$R^*(\overline\alpha): R^*(E^*) @>>> R^*(F^*)$ is a homotopy
equivalence between the minimal complexes of projectives.  Hence by
Lemma \ref{NAK} $R^*(\overline\alpha)$ is actually an isomorphism of
complexes.  So its naive truncation $\sigma_{\geq 0}
R^*(\overline\alpha) = \overline\alpha$ is also an isomorphism.
Therefore $\alpha$ is an isomorphism.
\end{pf}

We will also need a slight generalization of the previous lemma.

\begin{lemma}
\label{second:isom}
Let $\cal E = \cal H\zeta(\cal E)\oplus\bigoplus\cal O_{\PN}(n_i)$ and
$\cal F$ be vector bundles on $\PN$, and let $\cal Q$ be a coherent
sheaf on $\PN$.  Suppose that there exist morphisms $\alpha: \cal E
@>>> \cal F$ and $\beta: \cal F @>>> \cal Q$ such that

\rom(i\rom) $H^i_*(\alpha): H^i_*(\cal E) @>>> H^i_*(\cal F)$ is an
isomorphism for $0<i<N$,

\rom(ii\rom) $\beta\alpha$ takes the generators of the factors
$S(n_i)$ of $\Gamma_*(\cal E)$ onto a minimal set of generators of the
module $\overline Q :=\Gamma_*(\cal Q)/\beta\alpha(\Gamma_*(\cal
H\zeta(\cal E)))$,

\rom(iii\rom) $\cal E$ and $\cal F$ have the same rank.

Then $\alpha$ is an isomorphism.
\end{lemma}

\begin{pf}
Write $\cal F = \cal H\zeta(\cal F)\oplus\bigoplus\cal O_{\PN}(m_j)$.
The splittings of $\cal E$ and of $\cal F$ into direct factors are not
canonical.  But choosing such splittings gives an injection $\cal
H\zeta (\cal E) \hookrightarrow \cal E$ and a projection $\cal F
\twoheadrightarrow \cal H\zeta (\cal F)$.  Then the composition
\[
\overline{\alpha}:\quad \cal H\zeta (\cal E) @>>> \cal E @>{\alpha}>>
\cal F @>>> \cal H\zeta (\cal F)
\]
is, like $\alpha$, an isomorphism on $H^i_*$ for $0<i<N$.  So
$\overline \alpha$ is an isomorphism by Lemma \ref{isom}.  Hence by
identifying $\cal H\zeta(\cal F)$ with $\alpha (\cal H\zeta(\cal E))
\subset \cal F$, we see that $\alpha$ induces a morphism of diagrams
\begin{equation}
\label{split}
\begin{CD}
0 @>>> \cal H\zeta(\cal E) @>>> \cal E @>>> \bigoplus \cal
O_{\PN}(n_i) @>>> 0\\
&& @|  @VV{\alpha}V @VV{\alpha_1}V \\
0 @>>> \cal H\zeta(\cal F) @>>> \cal F @>>> \bigoplus \cal
O_{\PN}(m_j) @>>> 0
\end{CD}
\end{equation}
The morphisms $\alpha$ and $\beta$ therefore induce maps
\[
\bigoplus S(n_i) @>{\Gamma_*(\alpha_1)}>> \bigoplus S(m_j)
@>{\overline \beta}>> \overline Q = \Gamma_*(\cal Q)/ \beta\alpha
(\Gamma_*(\cal H\zeta(\cal E))).
\]
The composition is a surjection corresponding to a minimal set of
generators of $\overline Q$ by hypothesis (ii).  Hence the righthand
map $\overline \beta$ must be a surjection corresponding to a set of
generators of $\overline Q$.  However, the two free modules have the
same rank by hypothesis (iii).  Hence $\overline \beta$ also
corresponds to a minimal set of generators, and $\Gamma_*(\alpha_1)$
must be an isomorphism.  So returning to diagram (\ref{split}),
$\alpha_1$ and hence $\alpha$ are isomorphisms.
\end{pf}

\section{The Self-Dual Resolution}

Let $X\subset \Pnt$ be a locally Gorenstein subcanonical subscheme of
equicodimension $3$ satisfying the parity condition.  In this section
we use the Horrocks correspondence and especially Proposition
\ref{functorial} to construct a locally free resolution of $\cal O_X$.
We then use Lemma \ref{second:isom} to show that the resolution
satisfies condition (a) of Proposition
\ref{conditions}.

In the course of the construction we will need a more refined variant
of the canonical truncation.  Namely, suppose $D^*$ is a complex of
objects in an abelian category with differentials $\delta^i: D^i @>>>
D^{i+1}$.  Suppose $r$ is an integer, and $W\subset H^r(D^*)$ a
subobject.  Then $W$ may be pulled back to a $\overline W$ satisfying
\[
\im(\delta^{r-1})\subset \overline W \subset \ker(\delta^r) \subset
D^r
\]
We then define:
\[
\tau_{\leq r,W}(D^*): \qquad \dotsb @>>>  D^{r-2} @>>>  D^{r-1}
@>>> \overline{W} @>>> 0 @>>> 0 @>>> \dotsb \quad .
\]
The cohomology of this complex is given by
\[
H^i(\tau_{\leq r,W}(D^*))  =
\begin{cases}
H^i(D^*) & \text{if }i<r,\\ W & \text{if }i=r,\\0 & \text{if }i>r.
\end{cases}\\
\]

We will also use the following conventions.  If $\cal{E}$ is a
coherent sheaf on $\PN$ and $\alpha,\beta\in\Bbb{Q}$ are not both
integers, then we define $H^{\alpha}(\cal{E}(\beta)) = 0$.  Also if
$D^*$ is a complex and $\alpha\in\Bbb Q$, we define $\tau_{\leq
\alpha}(D^*) = \tau_{\leq [\alpha]}(D^*)$.

\subsection{Definition of the Locally Free Resolution}

Suppose $X\subset\Pnt$ is a locally Gorenstein subscheme of
equidimension $n>0$ such that $\omega_X\cong\cal O_X(l)$ for some
integer $l$ and such that $h^{n/2}(\cal O_X(l/2))$ is even.

Let $\nu = n/2$ and $l' = l/2$.  By hypothesis $H^{\nu}(\cal O_X(l'))$
is an even-dimensional vector space (zero if $n$ or $l$ is odd)
equipped with a nondegenerate $(-1)^{\nu}$-symmetric bilinear form
\[
H^{\nu}(\cal O_X(l')) \times H^{\nu}(\cal O_X(l')) @>>> H^n(\cal
O_X(l)) \cong k .
\]
Let $U\subset H^{\nu}(\cal O_X(l'))$ be an isotropic subspace of
maximal dimension $h^{\nu}(\cal O_X(l'))/2$.  Let
\begin{equation}
\label{w}
W = U \oplus \bigoplus_{t>l'} H^{\nu}(\cal O_X(t)) \subset
H^{\nu}_*(\cal O_X)
\end{equation}

We begin the construction of the locally free resolution with the
short exact sequence
\begin{equation}
\label{ideal}
0 @>>> \cal I_X @>>> \cal O_{\Pnt} @>>> \cal O_X @>>> 0.
\end{equation}
Since $H^i_*(\cal O_X)\cong H^{i+1}_*(\cal I_X)$ for $0<i<n+2$, we
have $W\subset H^{\nu+1}_*(\cal I_X)$.

Now since $X$ is locally Cohen-Macaulay of equidimension $n$, the
modules $H^i_*(\cal I_X)$ are of finite length for $0<i<n+1$.  Hence
the truncated complex $C^*_X = \tau_{>0} \tau_{\leq\nu+1,W}
\bold{R}\Gamma_*(\cal I_X)$ has cohomology modules $H^i(C^*_X)$ of
finite length for $0<i\leq \nu+1$, while $H^i(C^*_X)=0$ for all other
$i$.  Hence $C^*_X$ is in $\frak Z$.

The definition of $C^*_X$ as a truncation means that it is endowed
with a natural map $\beta: C^*_X @>>> \tau_{>0}\tau_{<n+3}
\bold{R}\Gamma_*(\cal I_X)$.  By Proposition \ref{functorial} this map
induces a morphism $\widetilde{\beta}: \cal{H}(C^*_X) @>>> \cal I_X$.
Let $Q$ be the cokernel
\[
H^0_*(\cal{H}(C^*_X)) @>{H^0_*(\widetilde{\beta})}>> H^0_*(\cal I_X)
@>>> Q @>>> 0.
\]
Let $d_1,\dots,d_r$ be the degrees of a minimal set of generators of
$Q$.  These generators lift to $H^0_*(\cal I_X)$, allowing us to
define a surjection
\begin{equation}
\label{gamma}
\gamma: \cal{F}_1 := \cal{H}(C^*_X)\oplus\bigoplus\cal O_{\Pnt}(-d_i)
\twoheadrightarrow \cal I_X .
\end{equation}
By construction, $\cal F_1$ is locally free.

Let $\cal K = \ker(\gamma)$.  We may then attach the short exact
sequence $0 @>>> \cal K @>>> \cal F_1 @>>> \cal I_X @>>> 0$ to the
short exact sequence (\ref{ideal}) to get an exact sequence
\begin{equation}
\label{partial}
0 @>>> \cal K @>>> \cal F_1 @>>> \cal O_{\Pnt} @>>> \cal O_X @>>> 0.
\end{equation}

The construction described above leads immediately to the following
conclusions about the cohomology of $\cal F_1$ and about the induced
morphisms $H^i_*(\gamma): H^i_*(\cal F_1) @>>> H^i_*(\cal I_X)$ (cf.\
Proposition \ref{functorial}).

\begin{itemize}
\item $H^i_*(\gamma)$ is surjective (resp.\ an isomorphism) for $i=0$
(resp.\ $0<i<\nu+1$).

\item $H^{\nu+1}_*(\gamma): H^{\nu+1}_*(\cal F_1) \cong W
\hookrightarrow H^{\nu+1}_*(\cal I_X)$ is injective.

\item $H^i_*(\cal F_1)=0$ for $\nu+1<i<n+3$.
\end{itemize}

One may now draw the following conclusions about the cohomology of
$\cal K$.

\begin{itemize}
\item $H^i_*(\cal K) = 0$ for $0<i<\nu+2$.

\item $H^{\nu+2}_*(\cal K)\cong H^{\nu}_*(\cal O_X)/W$.

\item $H^i_*(\cal K) \cong H^{i-2}_*(\cal O_X)$ for $\nu+2<i<n+3$.
\end{itemize}

To finish the definition of the locally free resolution, consider the
isomorphisms
\[
\Ext^1(\cal K,\omega_{\Pnt}(-l)) \cong H^{n+2}(\cal K(l))^* \cong
H^n(\cal O_X(l))^* \cong H^0(\cal O_X).
\]
The extension class corresponding to $1\in H^0(\cal O_X)$ gives a
short exact sequence
\begin{equation}
\label{extension}
0 @>>> \omega_{\Pnt}(-l) @>>> \cal F_2 @>>> \cal K @>>> 0
\end{equation}
which we may attach to (\ref{partial}) to get a complex of the type
(\ref{P}) resolving $\cal O_X$
\begin{equation}
\label{free:resol}
\cal P^*:\qquad 0 @>>> \omega_{\Pnt}(-l) @>>> \cal{F}_2 @>>> \cal{F}_1
@>>> \cal O_{\Pnt} .
\end{equation}

\begin{lemma}
\label{duality}
The sheaves $\cal F_1$ and $\cal F_2$ in the resolution
\rom{(\ref{free:resol})} satisfy $H^i_*(\cal F_2) \cong \left(
H^{n+3-i}_*(\cal F_1) \right)^*(l)$ for $0<i<n+3$.
\end{lemma}

\begin{pf}
If $0<i<\nu+2$, then $H^i_*(\cal F_2) \cong H^i_*(\cal K) = 0$ and
$H^{n+3-i}_*(\cal F_1) = 0$.  So the lemma holds for these values of
$i$.

If $i=\nu+2$, then $H^{\nu+2}_*(\cal F_2) \cong H^{\nu+2}_*(\cal K)
\cong H^{\nu}_*(\cal O_X)/W$, while $H^{\nu+1}_*(\cal F_1) \cong W$.
However, the submodule $W\subset H^{\nu}_*(\cal O_X)$ has been
constructed so that it is an isotropic submodule with respect to the
perfect pairing of Serre duality
\[
H^{\nu}_*(\cal O_X)\times H^{\nu}_*(\cal O_X) @>>> H^n_*(\cal
O_X) @>{\tr}>> k(-l).
\]
Moreover the length of $W$ is half the length of $H^{\nu}_*(\cal
O_X)$.  Hence $W = W^{\perp}$, and the duality isomorphism
$H^{\nu}_*(\cal O_X) \cong \left( H^{\nu}_*(\cal O_X) \right)^*(l)$
carries the submodule $W$ onto $\left(H^{\nu}_*(\cal O_X)/W
\right)^*(l)$.

If $\nu+2 < i < n+2$, then $H^i_*(\cal F_2) \cong H^i_*(\cal K) \cong
H^{i-2}_*(\cal O_X)$, while $H^{n+3-i}_*(\cal F_1) \cong
H^{n+3-i}_*(\cal I_X) \cong H^{n+2-i}_*(\cal O_X)$.  The asserted
duality is then simply the Serre duality pairing
\[
H^{i-2}_*(\cal O_X)\times H^{n+2-i}_*(\cal O_X) @>>> H^n_*(\cal
O_X) @>{\tr}>> k(-l).
\]

Finally if $i=n+2$, we have an exact sequence
\[
0 @>>> H^{n+2}_*(\cal F_2) @>>> H^{n+2}_*(\cal K) @>>>
H^{n+3}_*(\omega_{\Pnt}(-l)).
\]
Now $H^{n+2}_*(\cal K) \cong H^{n}_*(\cal O_X)$.  Moreover, the fact
that the extension class defining $\cal F_2$ corresponded under the
Serre duality identifications to $1\in H^0_*(\cal O_X) \cong \left(
H^{n+2}_*(\cal K) \right)^*(l)$ implies that the last exact sequence
dualizes to
\[
H^0_*(\cal O_{\Pnt}) @>1>> H^0_*(\cal O_X) @>>> \left( H^{n+2}_*(\cal
F_2) \right)^*(l) @>>> 0.
\]
Hence $\left( H^{n+2}_*(\cal F_2) \right)^*(l) \cong H^1_*(\cal I_X)
\cong H^1_*(\cal F_1)$.  Dualizing now gives the last of the asserted
isomorphisms.
\end{pf}

\begin{corollary}
The coherent sheaf $\cal F_2$ in the resolution
\rom{(\ref{free:resol})} is locally free.
\end{corollary}

\begin{pf}
Since $\cal F_1$ is locally free, $H^i_*(\cal F_1)$ is of finite
length for $0<i<n+3$.  So by the lemma, $H^i_*(\cal F_2)$ is also of
finite length for $0<i<n+3$.  But this implies that $\cal F_2$ is
locally free.
\end{pf}

\begin{proposition}
\label{cond:a}
The locally free resolution \rom{(\ref{free:resol})} satisfies
condition \rom(a\rom) of Proposition \ref{conditions}.
\end{proposition}

\begin{pf}
We write $\cal L=\omega_{\Pnt}(-l)$.  We have to show that if there is
a commutative diagram
\begin{equation}
\label{self:dual}
\begin{CD}
0 @>>> \cal{L} @>{d_3}>> \cal{F}_2 @>{d_2}>>
\cal{F}_1 @>{d_1}>> \cal O_{\Pnt}\\
&& @|  @VV{s_2}V  @VV{s_1}V  @| \\
0 @>>> \cal{L} @>{d_1\spcheck}>>
\cal{F}_1\spcheck\otimes\cal{L} @>{-d_2\spcheck}>>
\cal{F}_2\spcheck\otimes\cal{L} @>{d_3\spcheck}>> \cal O_{\Pnt}
\end{CD}
\end{equation}
such that the vertical maps extend the identity on $\cal O_X$, then
$s_1$ and $s_2$ are isomorphisms.

By exactness, the image of $d_3\spcheck$ is $\cal I_X$.  We will show
that $s_1$ is an isomorphism by applying Lemma \ref{second:isom} to
the composition
\[
\cal F_1 @>{s_1}>> \cal F_2\spcheck\otimes \cal L \twoheadrightarrow
\cal I_X.
\]
Note that this composition is exactly the surjection $\gamma: \cal F_2
\twoheadrightarrow \cal I_X$ of (\ref{gamma}).  Hence the composition
\[
H^i_*(\cal F_1) @>{H^i_*(s_1)}>> H^i_*(\cal F_2\spcheck\otimes\cal L) @>>>
H^i_*(\cal I_X)
\]
is injective for $0<i<n+3$.  A fortiori, $H^i_*(s_1)$ is also
injective for $0<i<n+3$.

However, by Serre duality $H^i_*(\cal F_2\spcheck\otimes \cal L)
\cong \left( H^{n+3-i}_*(\cal F_2) \right)^*(l)$ for all $i$.  So by
Lemma \ref{duality}, we have $H^i_*(\cal F_2\spcheck\otimes \cal L)
\cong H^i_*(\cal F_1)$ for $0<i<n+3$.  Hence for each $0<i<n+3$, the
morphism $H^i_*(s_1)$ is an injection of modules of the same finite
length.  Hence $H^i_*(s_1)$ is an isomorphism for $0<i<n+3$.  Thus
condition (i) of Lemma \ref{second:isom} holds.

Condition (ii) of Lemma \ref{second:isom} holds because of the method
of construction of $\cal F_1$ and of the surjection $\gamma$ in
(\ref{gamma}).  Finally exactness in the resolution implies that $\cal
F_1$ and $\cal F_2$ have the same rank.  Hence $\cal F_1$ and $\cal
F_2\spcheck\otimes\cal L$ also have the same rank, which is condition
(iii) of Lemma \ref{second:isom}.  Hence all three conditions of Lemma
\ref{second:isom} hold, and we may conclude that $s_1$ is an
isomorphism.

The map $s_2$ must now also be an isomorphism by the five-lemma.  This
completes the proof of the proposition.
\end{pf}

\section{The Differential Graded Algebra Structure}

In this section we finish the proof of Theorem \ref{main} by showing
that the locally free resolution (\ref{free:resol}) defined in the
previous section satisfies condition (b) of Proposition
\ref{conditions}.  That is to say, we show that the locally free
resolution (\ref{free:resol}) admits a commutative, associative
differential graded algebra structure.

Throughout this section we assume that the characteristic is not $2$.

We recall what needs to be proven.  In the previous section we defined
a locally free resolution (\ref{free:resol}) of $\cal O_X$
\[
0 @>>> \cal L @>{d_3}>> \cal F_2 @>{d_2}>> \cal F_1 @>{d_1}>> \cal
O_{\Pnt} @>>> \cal O_X.
\]
Let $\cal K=\ker(d_1)$.  We then had a morphism $\psi: \Lambda^2\cal
F_1 @>>> \cal K$ defined by $\psi(a\wedge b) = d_1(a)b-d_1(b)a$.  We
also have a long exact sequence
\[
\cdots @>>> \Hom(\Lambda^2\cal F_1,\cal F_2) @>>> \Hom(\Lambda^2\cal
F_1,\cal K) @>>> \Ext^1(\Lambda^2\cal F_1,\cal L) @>>> \cdots.
\]
According to diagram (\ref{liftdiag}), the problem is to lift $\psi\in
\Hom(\Lambda^2\cal F_1,\cal K)$ to a $\phi\in\Hom(\Lambda^2\cal
F_1,\cal F_2)$.  The obstruction to doing this is simply the image of
$\psi$ in
\[
\Ext^1(\Lambda^2\cal F_1,\cal L) = H^{n+2}(\Lambda^2\cal F_1(l))^*.
\]

Our first goal will therefore be to compute $H^{n+2}(\Lambda^2\cal
F_1(l))$.  We begin by considering a complex of locally free sheaves
on $\PN$.
\begin{equation}
\label{G}
\cal G^*: \qquad 0 @>>> \cal G^0 @>>> \cal G^1 @>>> \cdots @>>> \cal
G^r @>>> 0.
\end{equation}
There is an involution
\begin{align*}
T: \quad \cal G^*\otimes\cal G^* \ & @>>> \ \cal G^*\otimes\cal G^*\\
a\otimes b \ & \mapsto \ (-1)^{(\deg a)(\deg b)}b\otimes a
\end{align*}
interchanging the factors of $\cal G^*\otimes\cal G^*$.  Since the
characteristic is not $2$, the complex $\cal G^*\otimes\cal G^*$
splits into a direct sum of subcomplexes on which $T$ acts as
multiplication by $\pm 1$, viz.\ $\cal G^*\otimes\cal G^* = S_2(\cal
G^*) \oplus \Lambda^2(\cal G^*)$.  The complex $\Lambda^2(\cal G^*)$
is of the form
\begin{equation}
\label{Lambda:length}
\Lambda^2(\cal G^*):\qquad 0 @>>> \cal H^0 @>>> \cal H^1 @>>> \cdots
@>>> \cal H^{2r} @>>> 0
\end{equation}
where (cf.\ \cite{BE} p.\ 452)
\begin{equation}
\label{Lambda:formula}
\cal H^i \cong \bigoplus_{q<i/2} \left( \cal G^q\otimes \cal G^{i-q}
\right) \oplus
\begin{cases}
0 &\text{if $i$ is odd,}\\
\Lambda^2(\cal G^{i/2}) &\text{if }i\equiv 0\pmod 4,\\
S_2(\cal G^{i/2}) &\text{if }i\equiv 2\pmod 4.
\end{cases}
\end{equation}

\begin{lemma}
\label{Lambda}
Suppose $\cal G^*$ is a complex of locally free sheaves on $\PN$ as in
\rom{(\ref{G})} which is exact except in degree $0$.  Let $\cal E =
H^0(\cal G^*)$.  Then $\Lambda^2(\cal G^*)$ is an exact sequence of
locally free sheaves which is exact except in degree $0$, and
$H^0(\Lambda^2(\cal G^*)) = \Lambda^2\cal E$.
\end{lemma}

\begin{pf}
The standard spectral sequences of the double complex $\cal
G^*\otimes\cal G^*$ degenerate to show that the simple complex $\cal
G^*\otimes \cal G^*$ is exact except in degree $0$, and $H^0(\cal
G^*\otimes\cal G^*) = \cal E\otimes\cal E$.  Thus the augmented
complex $0 @>>> \cal E\otimes\cal E @>>> \cal G^*\otimes\cal G^*$ is
exact, and consequently its direct factor $0 @>>> \Lambda^2\cal E @>>>
\Lambda^2(\cal G^*)$ is also exact.
\end{pf}

\begin{lemma}
\label{max:cohom}
Suppose $\cal E$ is a locally free sheaf on $\PN$.  Let $r<N/2$ be an
integer.  Suppose that $H^i_*(\cal E)=0$ for $r<i<N$.  Then

\rom(a\rom) $H^i_*(\Lambda^2\cal E)=0$ for $2r<i<N$,

\rom(b\rom) $H^{2r}_*(\Lambda^2\cal E)\cong S_2 (H^r_*(\cal E))$ if
$r$ is odd, and $H^{2r}_*(\Lambda^2\cal E)\cong \Lambda^2 (H^r_*(\cal
E))$ if $r$ is even.

\rom(c\rom) If $H^r(\cal E(t))=0$ for $t<q$ for some integer $q$, then
$H^{2r}(\Lambda^2\cal E(t))=0$ for $t<2q$, while $H^{2r}(\Lambda^2\cal
E(2q)) \cong S_2 (H^r(\cal E(q)))$ if $r$ is odd, and $H^{2r}
(\Lambda^2 \cal E(2q)) \cong \Lambda^2 (H^r(\cal E(q)))$ if $r$ is
even.
\end{lemma}

\begin{pf}
By Lemma \ref{tauresol}(b), the minimal projective resolution of $P^*$
of the truncation $\tau_{<N}\bold R\Gamma_*(\cal E)$ is a complex of
free graded $S$-modules such that $P^i=0$ unless $0\leq i\leq N-1$.
Indeed, since $H^i(\tau_{<N}\bold R\Gamma_*(\cal E))=0$ for all $i>r$,
Lemma \ref{bounds} indicates that $P^i=0$ unless $0\leq i\leq r$,
i.e.\ $P^*$ is of the form
\[
P^*: \qquad 0 @>>> P^0 @>>> \cdots @>>> P^{r-1} @>>> P^r @>>> 0 .
\]

We now consider the complex of free graded $S$-modules
\[
\Lambda^2(P^*):\qquad 0 @>>> \Lambda^2 P^0 @>>> \cdots @>>> P^{r-1}
\otimes P^r @>>> T_2 (P^r) @>>> 0.
\]
where $T_2(P^r) = \Lambda^2(P^r)$ if $r$ is even, and $T_2(P^r) =
S_2(P^r)$ if $r$ is odd (cf.\ \eqref{Lambda:length} and
\eqref{Lambda:formula}).  According to Lemma \ref{tauresol}, the
complex of sheaves $\widetilde P^*$ associated to $P^*$ is exact
except in degree $0$ where the homology is $\cal E$.  So Lemma
\ref{Lambda} implies that the complex of sheaves $\Lambda^2(\widetilde
P^*)$ is also exact except in degree $0$ where the homology is
$\Lambda^2\cal E$.  The complex $\Lambda^2(P^*)$ of graded $S$-modules
therefore has homology of finite length except in degree $0$.
Moreover, the complex $\Lambda^2(P^*)$ vanishes except in degrees
between $0$ and $2r<N$, and the coefficients of its differentials lie
in $\frak m$ because it those of $P^*$ and therefore $P^*\otimes P^*$
do.  It now follows from Lemma \ref{tauresol}(a) that $\Lambda^2(P^*)$
is the minimal projective resolution of $\tau_{<N}\bold
R\Gamma_*(\Lambda^2 \cal E)$.

Therefore $H^i_*(\Lambda^2\cal E)\cong H^i(\Lambda^2(P^*))$ for all
$i<N$.  In particular, since $\Lambda^2(P^*)$ is concentrated in
degrees between $0$ and $2r$ by \eqref{Lambda:length}, we see that
$H^i_*(\Lambda^2\cal E)=0$ for $2r<i<N$.  This is part (a) of the
lemma.

For (b) note that $H^r_*(\cal E)$ and $H^{2r}_*(\Lambda^2\cal E)$ has
respective presentations
\[
\begin{CD}
P^{r-1} & \:\overset{\delta}{\longrightarrow}\: & P^r &
\:\rightarrow\: & H^r_*(\cal E) & \:\rightarrow\: & 0 ,\\
P^{r-1}\otimes P^r & \:\overset{\delta_1}{\longrightarrow}\: &
T_2(P^r) & \:\rightarrow\: & H^{2r}_*(\Lambda^2\cal E) &
\:\rightarrow\: & 0,
\end{CD}
\]
where $\delta_1(e\otimes f) = \delta(e)f\in T_2 (P^r)$.  But since the
presentation of $T_2(H^r_*(\cal E))$ is of exactly this form, we see
that $H^{2r}_*(\Lambda^2\cal E) \cong T_2(H^r_*(\cal E))$, as asserted
by the lemma.

For (c) write $H=H^r(\cal E(q))$.  The hypothesis that $H^r(\cal
E(t))=0$ for $t<q$ implies that
\(
P^r = \left( H\otimes_k S(-q)\right) \oplus F
\)
with $F=\bigoplus S(-n_i)$ for some $n_i>q$.  Then
\(
T_2(P^r) = \left( T_2 H \otimes_k S(-2q) \right) \oplus G
\)
with $G = \left( H\otimes_k F(-q) \right) \oplus T_2 F = \bigoplus
S(-m_j)$ for some $m_j>2q$.  Since the presentation of
$H^{2r}_*(\Lambda^2 \cal E)$ given above has the property that no
direct factor of $P^{r-1}\otimes P^r$ is mapped surjectively onto a
factor of $T_2(P^r)$, it now follows that $H^{2r}(\Lambda^2\cal
E(t))=0$ for $t<2q$, and $H^{2r}(\Lambda^2\cal E(2q)) \cong T_2 H$.
\end{pf}

\begin{corollary}
\label{Lambda:cohom}
Let $n$, $l$, and $X\subset \Pnt$ be as in Theorem \ref{main}.
Suppose that $U\subset H^{n/2}(\cal O_X(l/2))$ is the maximal
isotropic subspace defined in \eqref{w}, and that $\cal F_1$ is the
locally free sheaf defined in \eqref{gamma}.  Then
\[
H^{n+2}(\Lambda^2\cal F_1(l)) \cong \begin{cases} 0 &\text{if $n$ or
$l$ is odd,}\\ S_2 U &\text{if $l$ is even, and }n\equiv 0\pmod 4,\\
\Lambda^2 U &\text{if $l$ is even, and }n\equiv 2\pmod 4.
\end{cases}
\]
\end{corollary}

\begin{pf}
If $n$ is odd, then $H^i_*(\cal F_1)=0$ for $(n+1)/2 < i < n+3$.  So
Lemma \ref{max:cohom}(a) applies with $r=(n+1)/2$.  Therefore
$H^i_*(\Lambda^2\cal F_1) = 0$ for $n+1<i<n+3$, i.e.\
$H^{n+2}(\Lambda^2\cal F_1(t)) = 0$ for all $t$.

If $n$ is even but $l$ is odd, then Lemma \ref{max:cohom}(c) applies
with $r=(n+2)/2$ and $q=(l+1)/2$.  Then $H^{n+2}(\Lambda^2\cal
F_1(t))=0$ for all $t < l+1$.

If $l$ and $n$ are even, then Lemma \ref{max:cohom}(c) applies with
$r=(n+2)/2$ and $q=l/2$.  Since $H^{(n+2)/2}(\cal F_1(l/2))\cong U$,
it follows that $H^{n+2}(\Lambda^2\cal F_1(l))\cong \Lambda^2 U$ if
$r$ is even, and $H^{n+2}(\Lambda^2\cal F_1(l))\cong S_2 U$ if $r$ is
odd.  The corollary follows.
\end{pf}

\begin{lemma}
\label{cond:b}
If $\cal F_1$ is the locally free sheaf defined in \eqref{gamma}, then
the image of the map $\psi$ of \eqref{liftdiag} in
$\Ext^1(\Lambda^2\cal F_1,\cal L) \cong H^{n+2}(\Lambda^2\cal
F_1(l))^*$ vanishes.
\end{lemma}

\begin{pf}
If $n$ or $l$ is odd, then $H^{n+2}(\Lambda^2\cal F_1(l)) = 0$
according to Corollary \ref{Lambda:cohom}, so the image of $\psi$ is
evidently zero.

If $n$ and $l$ are even, then we claim that the image of $\psi$ in
$H^{n+1}(\Lambda^2\cal F_1(l))^*$ is the map
\[
\bigl\{ S_2 U \text{ or } \Lambda^2 U \bigr\} @>>> k
\]
which is the restriction to $U$ of the pairing $H^{n/2}(\cal O_X(l/2))
\times H^{n/2}(\cal O_X(l/2)) @>>> k$ of \eqref{pairing}.  Since $U$
was chosen isotropic, this map vanishes.

In order to prove the claim, we consider the diagonal $i: \Pnt =
\Delta\subset\Pnt\times\Pnt$.  Then there is a natural inclusion
$i(X)\subset X\times X$ which corresponds to a restriction map
\begin{equation}
\label{rest:diag}
\cal O_{X\times X} @>>> i_* \cal O_X.
\end{equation}
This map is essentially the multiplication $\cal O_X\otimes\cal O_X
@>>> \cal O_X$.  In any case applying $\bold R\Gamma_*$ to
\eqref{rest:diag} gives the cup product map
\begin{equation}
\label{cup:product}
\bold R\Gamma_*(\cal O_X)\otimes_k \bold R\Gamma_*(\cal O_X) @>>>
\bold R\Gamma_*(\cal O_X).
\end{equation}

Now consider the ``resolution'' of $\cal O_X$ given in \eqref{partial}
\[
\cal K^*:\qquad 0 @>>> \cal K @>>> \cal F_1 @>{d_1}>> \cal O_{\Pnt}
@>>> 0.
\]
The complex $\cal K^*$ is quasi-isomorphic to $\cal O_X$.  Hence the
restriction to the diagonal map \eqref{rest:diag} corresponds to a
morphism in the derived category
\[
p_1^*\cal K^*\otimes p_2^*\cal K^* @>>> i_*\cal K^*.
\]
In fact this morphism in the derived category is represented by an
actual map of complexes of sheaves
\begin{equation}
\label{FFK}
\begin{CD}
\dotsb & \:\rightarrow\: & p_1^*\cal K \oplus (p_1^*\cal F_1 \otimes
p_2^*\cal F_1) \oplus p_2^*\cal K & \:\rightarrow\: & p_1^*\cal F_1
\oplus p_2^*\cal F_1 & \:\rightarrow\: & \cal O_{\Bbb P\times\Bbb P} &
\:\rightarrow\: & 0\\
&& @VVV @VVV @VVV \\
0 & \:\rightarrow\: & i_*\cal K & \:\rightarrow\: & i_*\cal F_1 &
\:\rightarrow\: & i_*\cal O_{\Bbb P} & \:\rightarrow\: & 0
\end{CD}
\end{equation}
All the vertical maps are straightforward restrictions to the diagonal
except for the component $p_1^*\cal F_1 \otimes p_2^*\cal F_1 @>>>
i_*(\cal K)$ which is defined (like $\psi$ of
\eqref{liftdiag}) by noting that the composition
\begin{equation}
\label{FFK2}
\begin{CD}
p_1^*\cal F_1 \otimes p_2^*\cal F_1 & \:\rightarrow\: p_1^*\cal F_1
\oplus p_2^*\cal F_1 \:\rightarrow\: & i_* \cal F_1 \\
p_1^*(a)\otimes p_2^*(b) & \mapsto & i_*\left(d_1(a)b-d_1(b)a\right)
\end{CD}
\end{equation}
is contained in the kernel of $i_*\cal F_1 @>>> i_*\cal O_{\Bbb P}$.

Now since $\cal K^*$ is quasi-isomorphic to $\cal O_X$, if we apply
$\bold R\Gamma_*$ to \eqref{FFK} we get a morphism $\bold
R\Gamma_*(p_1^*\cal K\otimes p_2^*\cal K) @>>> \bold R\Gamma_*(i_*\cal
K^*)$ in $D^b_{\grssmod}$ which is quasi-isomorphic to
\eqref{cup:product}.  In particular, the maps of hypercohomology are
quasi-isomorphic to the cup product
\[
H^n_*(p_1^*\cal O_X\otimes p_2^*\cal O_X) \cong \bigoplus_i H^i_*(\cal
O_X)\otimes_k H^{n-i}_*(\cal O_X) @>>> H^n_*(\cal O_X).
\]

The hypercohomology $H^n_*(\cal K^*)\cong H^n_*(\cal O_X)$ is of
course the same as the $H^n$ of the total complex of the double
complex
\[
0 @>>> \bold R\Gamma_*(\cal K) @>>> \bold R\Gamma_*(\cal F_1) @>>>
\bold R\Gamma_*(\cal O_{\Bbb P}) @>>> 0.
\]
According to the calculations at the beginning of the previous
section, this $H^n_*$ is all attributable to $\cal K$, i.e.\ the
truncation $\cal K^* @>>> \cal K[2]$ induces an isomorphism
$H^n_*(\cal O_X) \cong H^n_*(\cal K^*) \cong H^{n+2}_*(\cal K)$.

Similarly, the hypercohomology $H^n_*(p_1^*\cal K^*\otimes p_2^*\cal
K^*)$ is the same as the $H^n$ of the total complex of $\bold
R\Gamma_*$ of the first row of \eqref{FFK}.  The submodule $W\otimes_k W
\subset H^n_*(p_1^*\cal O_X\otimes p_2^*\cal O_X)$ is
attributable as the $H^{n+2}_*$ of the factor $p_1^*\cal F_1\otimes
p_2^*\cal F_1$ in the first row of \eqref{FFK}.  Therefore $H^{n+2}_*$
of the vertical map $p_1^*\cal F_1 \otimes p_2^*\cal F_1 @>>> \cal
K$ is simply the cup product map $W\otimes_k W @>>> H^n_*(\cal O_X)$.

Now the fact that $i_*(\cal K)$ is supported on $\Delta$, plus the
symmetry of the product map imply that the vertical map of \eqref{FFK}
factors as
\begin{equation}
\label{factorization}
\begin{CD}
p_1^*\cal F_1\otimes p_2^*\cal F_1 & \:\rightarrow\: & i_*(\cal
F_1\otimes \cal F_1) & \:\rightarrow\: & i_*(\Lambda^2\cal F_1)
& \:\overset{i_*(\psi)}{\longrightarrow}\: & i_*(\cal K).\\
p_1^*(a)\otimes p_2^*(b) & \mapsto & i_*(a\otimes b) & \mapsto &
i_*(a\wedge b) & \mapsto & d_1(a)b-d_1(b)a
\end{CD}
\end{equation}
We wish to calculate $H^{n+2}_*$ of the above morphisms.  Let
\[
P^*:\qquad 0 @>>> P^0 @>>> \dotsb @>>> P^{(n+2)/2} @>>> 0
\]
be a minimal projective resolution of $\tau_{<n+3}\bold R\Gamma_*(\cal
F_1)$ (cf.\ Lemma \ref{tauresol}).  Then if one applies
$\tau_{<n+3}\bold R\Gamma_*$ to the first two morphisms of
\eqref{factorization}, one gets the natural maps
\[
P^*\otimes_k P^* @>>> P^*\otimes_S P^* @>>> \Lambda^2(P^*)
\]
(cf.\ the proof of Lemma \ref{max:cohom}).  All three complexes are
supported in degrees between $0$ and $n+2$, and applying $H^{n+2}$
gives surjections
\[
W\otimes_k W \twoheadrightarrow W\otimes_S W \twoheadrightarrow
\bigl\{ S_2 W \text{ or } \Lambda^2 W \bigr\}.
\]
It therefore follows that $H^{n+2}_*(\psi): H^{n+2}_*(\Lambda^2\cal
F_1) @>>> H^{n+2}_*(\cal K)$ is isomorphic to the cup product map
$\bigl\{ S_2 U \text{ or } \Lambda^2 U \bigr\} @>>> H^n_*(\cal O_X)$.
In particular, in degree $l$ the morphism $H^{n+2}(\cal F_1(l)) @>>>
H^{n+2}(\cal K(l))$ is the same as $\bigl\{ S_2 U \text{ or }
\Lambda^2 U \bigr\} @>>> H^n(\cal O_X(l))$.

We now have to consider the extension of \eqref{extension}
\[
0 @>>> \omega_{\Bbb P}(-l) @>>> \cal F_2 @>>> \cal K @>>> 0.
\]
(Recall $\cal L = \omega_{\Bbb P}(-l)$.)  In the associated long
exact sequence of cohomology
\[
\dotsb @>>> H^{n+2}(\cal F_2(l)) @>>> H^{n+2}(\cal K(l)) @>\tr>>
H^{n+3}(\omega_{\Pnt}) \cong k,
\]
the differential is the element of $H^{n+2}(\cal K(l))^* \cong
\Ext^1(\cal K(l),\omega_{\Pnt})$ corresponding to the extension class.
So by construction the differential is the trace map $\tr\in
H^{n+2}(\cal K(l))^*\cong H^n(\cal O_X(l))^*$ which corresponds under
Serre duality to $1\in H^0(\cal O_X)$.

Now the image of $\psi\in \Ext^1(\Lambda^2\cal F_1(l),\omega_{\Bbb P})
\cong H^{n+2}(\Lambda^2\cal F_1(l))^*$ is exactly the composition
\[
H^{n+2}(\Lambda^2\cal F_1(l)) @>{H^{n+2}(\psi)}>> H^{n+2}(\cal K(l))
@>\tr>> H^{n+3}(\omega_{\Pnt}) \cong k.
\]
By our previous calculations, this is the composition of the cup
product map $\bigl\{ S_2 U \text{ or } \Lambda^2 U \bigr\} @>>>
H^n(\cal O_X(l))$ with the trace map $H^n(\cal O_X(l)) @>>> k$.
Therefore this composition is the restriction to $S_2 U$ or $\Lambda^2
U$ of the Serre duality pairing $H^{n/2}(\cal O_X(l/2))\otimes
H^{n/2}(\cal O_X(l/2)) @>>> k$.  This is what was claimed at the
beginning of the proof of the lemma.  Since $U$ was chosen isotropic,
this composition vanishes, i.e.\ the image of $\psi$ in $\Ext^1(\cal
K,\cal L)$ vanishes.
\end{pf}

\begin{pf*}{Proof of Theorem \ref{main}}
According to Proposition \ref{conditions}, in order to prove Theorem
\ref{main} it suffices to find a locally free resolution
\[
0 @>>> \cal L @>>> \cal F_2 @>>> \cal F_1 @>>> \cal O_{\Pnt} @>>> \cal O_X
\]
which satisfies two conditions.  But the locally free resolution
defined in \eqref{free:resol} was shown to satisfy the first of these
conditions was shown in Proposition \ref{cond:a}.  Moreover, this
resolution was just shown to satisfy the second condition in Lemma
\ref{cond:b}.  Hence Theorem \ref{main} holds.
\end{pf*}

\section{Characteristic $2$ Computations}

In the introduction, we asserted that Theorem \ref{main} also holds in
characteristic $2$ provided the phrase ``$n\equiv 0\pmod 4$'' in the
parity condition is replaced by the phrase ``$n$ is even.''  In this
section we justify that assertion by proving analogues of Lemmas
\ref{Lambda} and \ref{max:cohom} and Corollary \ref{Lambda:cohom} in
characteristic $2$.  These were the only steps in the proof of
Theorem \ref{main} where we used the assumption that the
characteristic is not $2$.

Throughout this section we assume that the characteristic is $2$.

We recall certain simple facts from modular representation theory.
Let $R$ be a commutative algebra over a field of characteristic $2$,
and let $V$ be a free $R$-module.  Let $t\in\End(V\otimes V)$ be the
endomorphism $t(a\otimes b) = a\otimes b - b\otimes a$.  Set $D_2 V
=\ker(t)$, and $\Lambda^2 V =\im(t)$, and $S_2 V = \coker(t)$.  Since
$t^2=0$ in characteristic $2$, there are inclusions
\[
0 \subset \Lambda^2 V \subset D_2 V \subset V\otimes V
\]
and corresponding surjection of quotients of $V\otimes V$
\[
V\otimes V \twoheadrightarrow S_2 V \twoheadrightarrow \Lambda^2 V
@>>> 0.
\]
The subquotient $D_2 V /\Lambda^2 V$ is $F(V)$, the Frobenius pullback
of $V$.  It is a free module of the same rank as $V$.  This $F(V)$ is
also the kernel of the surjection $S_2 V \twoheadrightarrow \Lambda^2
V$.

Note that the natural map from $S_2 V $ to $V\otimes V$ given by $xy
\mapsto x\otimes y + y\otimes x$ is not injective because any $x^2
\mapsto 0$.  The map is the composition $S_2 V  \twoheadrightarrow
\Lambda^2 V \hookrightarrow V\otimes V$ with kernel $F(V)$.

The operations $D_2$, $\Lambda^2$, $S_2$, and $F$ are all functorial.
Therefore we may define $D_2\cal E$, $\Lambda^2\cal E$, $S_2\cal E$,
and $F(\cal E)$ for any locally free sheaf $\cal E$ on any scheme $X$
over a field of characteristic $2$.

In some ways $F$ has better properties than the others.  If $M =
(m_{ij}): R^n @>>> R^m$ is a morphism of free $R$-modules, then $F(M)
= (m_{ij}^2): R^n @>>> R^m$.  One may use this formula together with
the Buchsbaum-Eisenbud exactness criterion \cite{BE:exact} to show
that $F$ of an exact sequence of locally free sheaves is exact.  As
a result of this we get the following lemma.

\begin{lemma}
\label{frob}
Let $\cal E$ be a locally free sheaf on $\PN$.  If for some integer
$r$ one has $H^i_*(\cal E)=0$ for $r<i<N$, then $H^i_*(F(\cal E))=0$
for $r<i<N$ also.
\end{lemma}

\begin{pf}
Let $P^*$ be the minimal projective resolution of $\tau_{<N}\bold
R\Gamma_*(\cal E)$.  By Lemmas \ref{tauresol} and \ref{bounds}, $P^*$
has the form
\[
P^*:\qquad 0 @>>> P^0 @>>> P^1 @>>> \dotsb @>>> P^r @>>> 0.
\]
Moreover, $P^*$ is exact except in degree $0$ away from the irrelevant
ideal $\frak m\subset S$, and $H^0(\widetilde P^*)=\cal E$.

The functoriality and exactness of $F$ now imply that
\[
F(P^*):\qquad 0 @>>> F(P^0) @>>> F(P^1) @>>> \dotsb @>>> F(P^r) @>>> 0
\]
is exact except in degree in degree $0$ away from the irrelevant
ideal, and has $H^0(F(\widetilde P^*))=F(\cal E)$.  Applying Lemma
\ref{tauresol} again, we conclude that $F(P^*)$ is the minimal projective
resolution of $\tau_{<N}\bold R\Gamma_*(F(\cal E))$.

So if $r<i<N$, then $H^i_*(F(\cal E)) = H^i(F(P^*))=0$ since $F(P^*)$
vanishes in degrees greater than $r$.
\end{pf}

\begin{corollary}
\label{l:s}
Let $\cal E$ be a locally free sheaf on $\PN$ over a field of
characteristic $2$ such that $H^{N-1}_*(\Lambda^2\cal E)=0$.
Suppose $r<N$ is an integer such that $H^i_*(\cal E)=0$ for $r<i<N$.
Then $H^i_*(S_2\cal E) \cong H^i_*(\Lambda^2\cal E)$ for $r<i<N$.
\end{corollary}

\begin{pf}
We consider the exact sequence $0 @>>> F(\cal E) @>>> S_2\cal E @>>>
\Lambda^2\cal E @>>> 0$ and the associated long exact sequence
\[
\dotsb @>>> H^i_*(F(\cal E)) @>>> H^i_*(S_2\cal E) @>>>
H^i_*(\Lambda^2\cal E) @>>> H^{i+1}_*(F(\cal E)) @>>> \dotsb
\]
The hypothesis $H^i_*(\cal E)=0$ for $r<i<N$ implies also
$H^i_*(F(\cal E))=0$ for $r<i<N$ by Lemma \ref{frob}.  Hence the long
exact sequence implies that $H^i_*(S_2\cal E)\cong
H^i_*(\Lambda^2\cal E)$ for $r<i<N-1$ and that $H^{N-1}_*(S_2\cal
E) \hookrightarrow H^{N-1}_*(\Lambda^2\cal E)$ is injective.  But
by hypothesis $H^{N-1}_*(\Lambda^2\cal E)=0$, so $H^{N-1}_*(S_2\cal
E)=0$ as well.
\end{pf}

We now prove the analogue of Lemma \ref{max:cohom}.

\begin{lemma}
\label{char:two}
Suppose $\cal E$ is a locally free sheaf on $\PN$ over a field of
characteristic $2$.  Let $0 < r <N/2$ be an integer such that
$H^i_*(\cal E)=0$ for $r<i<N$.  Then

\rom(a\rom) $H^i_*(\Lambda^2\cal E)=0$ for $2r<i<N$,

\rom(b\rom) $H^{2r}_*(\Lambda^2\cal E)\cong S_2 (H^r_*(\cal E))$.

\rom(c\rom) If $H^r(\cal E(t))=0$ for $t<q$ for some integer $q$, then
$H^{2r}(\Lambda^2\cal E(t))=0$ for $t<2q$, while $H^{2r}(\Lambda^2\cal
E(2q)) \cong S_2 (H^r(\cal E(q)))$.
\end{lemma}

\begin{pf}
Let $P^*$ be the minimal projective resolution of $\tau_{<N}\bold
R\Gamma_*(\cal E)$
\[
P^*:\qquad 0 @>>> P^0 @>{\delta^0}>> P^1 @>{\delta^1}>> P^2 @>>>
\dotsb @>>> P^r @>>> 0
\]
Let $\cal P = \widetilde P^0$, and let $\cal F = \widetilde{\ker(
\delta^1)}$.  Then we have an exact sequence $0 @>>> \cal E @>>> \cal
P @>>> \cal F @>>> 0$ such that $\cal P$ is a direct sum of line
bundles, $H^i_*(\cal F) = H^{i+1}_*(\cal E)$ for $0<i<N-1$, and
$H^{N-1}_*(\cal F)=0$.

It is easy to see that there is a natural exact complex
\begin{equation}
\label{lambda:s}
0 @>>> \Lambda^2 \cal E @>>> \Lambda^2 \cal P @>>> \cal P\otimes
\cal F @>>> S_2\cal F @>>> 0.
\end{equation}

We now prove parts (a) and (b) of the lemma by induction on $r$.  If
$r=1$, then $\cal F = \widetilde P^1$ is a direct sum of line bundles,
and the complex \eqref{lambda:s} is just the augmented complex
\[
0 @>>> \Lambda^2 \cal E @>>> \Lambda^2(\widetilde P^*)
\]
which is still exact in this case.  So we may conclude just as in
Lemma \ref{max:cohom} that $H^i_*(\Lambda^2\cal E)=0$ for $2<i<N$,
and that $H^2_*(\Lambda^2\cal E) = S_2(H^1_*(\cal E))$.

If $r>1$, then $H^i_*(\cal F)=0$ for $r-1<i<N$.  So by induction
$H^i_*(\Lambda^2\cal F)=0$ for $2r-2<i<N$, and also $H^{2r-2}_*
(\Lambda^2\cal F)) \cong S_2(H^{r-1}_*(\cal F)) \cong S_2(H^r_*(\cal
E))$.  It now follows from Corollary \ref{l:s} that $H^i_*(S_2\cal
F)\cong H^i_*(\Lambda^2\cal F)=0$ for $r-1 < i < N$.  So in particular
$H^i_*(S_2\cal F)=0$ for $2r-2<i<N$ and that $H^{2r-2}_*(S_2\cal
F)\cong S_2(H^r_*(\cal E))$.  Now since $\cal P$ is a direct sum of
line bundles, we have $H^i_*(\Lambda^2\cal P)=0$ for $0<i<N$, and
$H^i_*(\cal P\otimes\cal F)=0$ for $r-1<i<N$.  So if we break up
\eqref{lambda:s} into short exact sequences and take its graded
cohomology, we can deduce that $H^i_*(\Lambda^2\cal E)) \cong
H^{i-2}_* (S_2\cal F)$ for $r+1<i<N$.  Since $r+1<2r$, this gives
parts (a) and (b) of the lemma.

Part (c) of the lemma follows from part (b) by the same argument as in
Lemma \ref{max:cohom}.
\end{pf}

We have the following corollary in analogy with Corollary
\ref{Lambda:cohom}.

\begin{corollary}
\label{S2U}
Let $n$, $l$, and $X\subset \Pnt_k$ be as in Theorem \ref{main} with
$k$ a field of characteristic $2$.  Suppose that $U\subset
H^{n/2}(\cal O_X(l/2))$ is the maximal isotropic subspace defined in
\eqref{w}, and that $\cal F_1$ is the locally free sheaf defined in
\eqref{gamma}.  Then
\[
H^{n+2}(\Lambda^2\cal F_1(l)) \cong \begin{cases} 0 &\text{if $n$ or
$l$ is odd,}\\ S_2 U &\text{if $n$ and $l$ are even.} \end{cases}
\]
\end{corollary}

The proof of Lemma \ref{cond:b} in characteristic $2$ is essentially
the same as in the previous section, only with Lemma \ref{char:two}
and Corollary \ref{S2U} replacing their analogues, and with $T_2 =
S_2$ always.  Hence Theorem \ref{main} also holds in characteristic
$2$ as long as one treats all even $n$ the same.

\section{The Local Version of the Main Theorem}

In this section we consider Theorem \ref{RLR}, the local version of
our main result.  We state a variant version which is clearly a local
analogue of Theorem \ref{main} with an identical proof, and then show
that this variant version is equivalent to Theorem \ref{RLR}.

Let $(R,\frak m,k)$ be a regular local ring, and let $U = \Spec(R) -
\{\frak m\}$ be the punctured spectrum of $R$.  We say that a closed
subscheme $Y\subset U$ of pure codimension $3$ is {\em Pfaffian} if
$O_X$ has a locally free resolution on $U$
\[
0 @>>> \cal{O}_U @>h>> \cal{E}\spcheck @>f>> \cal{E} @>g>> \cal{O}_U
@>>> \cal O_X
\]
where $\cal E$ is a locally free $\cal O_U$-module of odd rank $2p+1$,
$f$ is skew-symmetric, and $g$ and $h=g\spcheck$ are given locally by
the Pfaffians of order $2p$ of $f$.  The following theorem is the
obvious local analogue of Theorem \ref{main}.

\begin{theorem}
\label{punc:spec}
Let $(R,\frak m,k)$ be a regular local ring of dimension $n+4>4$ with
residue field not of characteristic $2$, and let $U=\Spec(R)-\{\frak
m\}$.  Let $X\subset U$ be a closed subscheme of pure codimension $3$.
Then $X$ is Pfaffian if and only if the following three conditions
hold:

\rom(a\rom) $X$ is locally Gorenstein,

\rom(b\rom) $\omega_X\cong \cal O_X$, and

\rom(c\rom) if $n\equiv 0\pmod 4$, then $H^{n/2}(\cal O_X)$ is of even
length.
\end{theorem}

Theorem \ref{punc:spec} may be proven in exactly the same manner as
Theorem \ref{main}.  All results concerning the Buchsbaum-Eisenbud
proof, the Horrocks correspondence, Serre/local duality, the
cohomology of $H^{n+2}(\Lambda^2\cal F_1)$ work identically for graded
modules over polynomial rings over $k$ and for modules over regular
local $k$-algebras.  There is only one point which is in any way more
subtle in the local case.  Namely, if $n$ is even, then one has a
Matlis duality pairing of $R$-modules of even finite length
\[
H^{n/2}(\cal O_X) \times H^{n/2}(\cal O_X) @>>> k.
\]
This pairing is perfect in the sense that for any submodule $M\subset
H^{n/2}(\cal O_X)$ one has
\[
\length(M)+\length(M^\perp) = \length(H^{n/2}(\cal O_X)).
\]
In order to be able to define $C^*_X$ and $\cal F_1$ as in
\eqref{gamma} one must choose an isotropic submodule $W$ of
length equal to half that of $H^{n/2}(\cal O_X)$.  But it is not
difficult to show that this is possible.

We now compare Theorems \ref{RLR} and \ref{punc:spec}.  First of all,
$E=\Gamma(\cal E)$ gives a bijective correspondence between locally
free sheaves $\cal E$ on $U$ and reflexive $R$-modules $E$ such that
$E_{\frak p}$ is a free $R_{\frak p}$-module for all prime ideals
$\frak p\neq \frak m$.  There is also bijective correspondence
betweenclosed subschemes $X\subset U$ of pure codimension $3$ and
unmixed ideals $I\subset R$ of height $3$ given by $I=\Gamma(\cal
I_X)$.  Hence an ideal $I$ is Pfaffian in the sense of Theorem
\ref{RLR} if and only if the corresponding subscheme $X\subset U$ is
Pfaffian in the sense of Theorem \ref{punc:spec}.

The three conditions (a), (b), and (c) of the two theorems also
correspond.  In the case of (a) this is obvious.  For (b) note that
$\omega_{R/I}\cong \Gamma(\omega_X)$ since for all $\frak p\in U$ one
has $\omega_{R/I,\frak p} = \Ext^3_{R_{\frak p}}((R/I)_{\frak
p},R_{\frak p}) = \omega_{X,\frak p}$, and $\omega_{R/I}$ is
saturated.  Similarly $(R/I)^{\sat} \cong \Gamma(\cal O_X)$.  This
gives the equivalence of the two conditions (b).

As for the conditions (c), first note that the dimension $n$ in
Theorem \ref{RLR} corresponds to $n+4$ in Theorem \ref{punc:spec}.
But if one uses $n$ as in the latter theorem, one has
\[
H^{n/2}(U,\cal O_X) \cong H^{(n+2)/2}(U,\cal I_X) \cong
H^{(n+4)/2}_{\frak m}(I).
\]
Hence the two conditions (c) correspond.

Therefore the two theorems \ref{RLR} and \ref{punc:spec} are
equivalent, as claimed.

In equicharacteristic $2$ the computations of Section 5 remain true in
the local case.  So Theorems \ref{RLR} and \ref{punc:spec} are true in
equicharacteristic $2$ provided one changes the phrase ``$n\equiv
0\pmod 4$'' in the parity condition to ``$n$ is even.''  If $R$ is a
regular local ring with residue field of characteristic $2$ and
quotient field of characteristic $0$, a different set of calculations
is needed.  These are unfortunately somewhat involved, and we do not
reproduce them here.


\begin{thebibliography}{MM9}
\bibitem[BE1]{BE:exact} D.~Buchsbaum and D.~Eisenbud, {\em What Makes
a Complex Exact}?, J.~Algebra {\bf 25} (1973), 259--268.

\bibitem[BE2]{BE} \bysame, {\em Algebra Structures for Finite Free
Resolutions, and Some Structure Theorems for Ideals of Codimension
$3$}, Amer.\ J.\ Math.\ {\bf 99} (1977), 447--485.

\bibitem[Ha]{Ha} R.~Hartshorne, {\em Residues and Duality}, Lecture
Notes in Math., vol.~20, Springer-Verlag, Berlin, 1966.

\bibitem[Ho]{H} G.~Horrocks, {\em Vector Bundles on the Punctured
Spectrum of a Local Ring}, Proc.\ London Math.\ Soc.\ (3) {\bf 14}
(1964), 689--713.

\bibitem[O]{O} C.~Okonek, {\em Notes on Varieties of Codimension $3$ in
$\PN$}, preprint, Bonn, 1993.

\bibitem[W]{W} C.~Walter, {\em Algebraic Methods for the Cohomology of
Normal Bundles of Algebraic Space Curves}, preprint.
\end{thebibliography}
\end{document}